\documentclass[12pt,leqno]{article}

\newcounter{conjecture}
\newcounter{remark}
\newcounter{corollary}

\newenvironment{corollary}{\medskip{\bf Corollary \thecorollary.}
\addtocounter{corollary}{1}\em}{\rm}
\newcommand{\eqnsection}{
     \renewcommand{\theequation}{\thesection.\arabic{equation}}
     \makeatletter
     \csname @addtoreset\endcsname{equation}{section}
     \makeatother}
\newtheorem{theorem}{Theorem}
\newtheorem{lemma}{Lemma}
\newtheorem{proposition}{Proposition}

\newcommand{\upa}{\uparrow}
\newcommand{\doa}{\downarrow}
\newcommand{\dd}{\delta}
\newcommand{\subs}{\subseteq}
\newcommand{\lar}{\longrightarrow}
\newcommand{\eps}{\varepsilon}

\newcommand{\aaa}{\alpha}
\newcommand{\reals}{R}
\newcommand{\inn}{\bigcap}
\newcommand{\unn}{\bigcup}
\newcommand{\lll}{\label}

\def \be{\begin{equation}}
\def \ee{\end{equation}}
\def \bt{\begin{theorem}}
\def \et{\end{theorem}}
\def \bea{\begin{eqnarray}}
\def \eea{\end{eqnarray}}
\def \bas{\begin{eqnarray*}}
\def \eas{\end{eqnarray*}}

\def \al{\alpha}
\def \bb{\beta}

\def \th{\theta}

\def \({\left(}
\def \){\right)}

\def \nn{\nonumber}

\def \bc{\begin{center} }
\def \ec{\end{center} }
\def \bs{\begin{slide} }
\def \es{\end{slide} }

\def\square{{\vcenter{\vbox{\hrule height.3pt
          \hbox{\vrule width.3pt height5pt \kern5pt
             \vrule width.3pt}
          \hrule height.3pt}}}}
\def\qed{{\hfill $\square$ \bigskip}}

\eqnsection
\begin{document}

\title{Renormalization and convergence in law for the derivative of intersection local time in ${\bf R}^2$}

\author{Greg Markowsky}

\bibliographystyle{amsplain}
\maketitle

\begin{abstract}
In this paper we will examine the derivative of intersection local
time of Brownian motion and symmetric stable processes in $R^2$.
These processes do not exist when defined in the canonical way. The
purpose of this paper is to exhibit the correct rate for
renormaliztion of these processes.

\end{abstract}

\eqnsection
\setlength{\unitlength}{2mm}

\section{Introduction}
Let $B_t$ be a Brownian motion in {\bf R}$^d$, and let

\be \lll{nopc} \al_{\eps}(T) =  \int_0^T \int_0^t
f_{\eps}(B_s - B_u) ds dt
\ee where $f_{\eps}$ denotes the Gaussian density function on ${\bf R}^{d}$ with variance $\eps$. If $\al_{\eps}(T)$ converges to a process as $\eps \lar 0$ we denote this process $\al_T$ and
call it the {\it intersection local time}(henceforth abbreviated as ILT).

In one dimension the ILT does exist, as can be seen easily by using
the occupation times formula(see, e.g. \cite{rosen2}). In dimension
$2$ it does not exist as defined above, as $\al_{\eps}(T)$ blows up
as $\eps \lar 0$ due to the set $\{s = u \}$. In \cite{varad},
however, Varadhan showed that $\al_{\eps}(T) - E[\al_{\eps}(T)]$
does converge in law to a process, which is referred to as the {\it
renormalized intersection local time}. This process was originally
considered due to its relevance to quantum field theory(see
\cite{varad}), but it has found several other uses, for example in
Edwards' work on polymers (see \cite{edwards}) and Le Gall's work on
Wiener sausages (see \cite{legall1}, \cite{legall2}).

In \cite{yor}, Yor proved that in d = 3

\be \Big{ \{ }\frac{1}{\sqrt{\log(1/\eps)}}\Big( \al_{\eps}(T) -
E[\al_{\eps}(T)] \Big), T \geq 0 \Big{\}} \ee converges in law as
$\eps \lar 0$ to the process $\{ \frac{1}{\sqrt{2}\pi} B_T, T \geq 0
\}$, where $B_T$ is a one dimensional Brownian motion. This theorem
inspired a similar result from Rosen in \cite{rosen1} involving
symmetric stable processes. Rosen considered the process

\be \aaa_{\eps}(T)= \int_{0}^{T} \int_{0}^{t}
f_{\eps}(X_{t}-X_{s})dsdt \lll{1}\ee where now $X$ is a symmetric
stable process of index $\bb$, and $f_{t}$ denotes the density of $X_t$. Rosen proved that if $4/3 < \bb \leq 2$, then $\al_{\eps}(T) - E[\al_{\eps}(T)]$ converges pathwise
as $\eps \lar 0$ to a finite random variable. If $\bb = 4/3$, then

\be \Big{ \{ } \frac{1}{\sqrt{\log(1/\eps)}}\Big( \al_{\eps}(T) -
E[\al_{\eps}(T)] \Big), T \geq 0 \Big{ \} } \ee converges in law as
$\eps \lar 0$ to $\{k(\bb) B_T, T \geq 0\}$ where $k(\bb)$ is a
constant which depends on $\bb$. Similarly, if $1< \bb < 4/3$, then

\be \Big{ \{ }\eps^{2/\bb - 3/2} \Big( \al_{\eps}(T) -
E[\al_{\eps}(T)] \Big), T \geq 0 \Big{ \} }\ee converges in law as
$\eps \lar 0$ to $\{k(\bb) B_T, T \geq 0\}$ where $k(\bb)$ is a
constant depending on $\bb$. When $\bb = 2$, $X$ is Brownian motion, and this gives a different
proof of Varahdan's renormalization. The method employed by Rosen in proving this
also gives an alternate proof of Yor's result in 3 dimensions.

In \cite{rosen2} Rosen introduced the notion of the
derivative of the intersection local time of Brownian motion in
${\bf R}^1$. It is defined as

\be \lll{nopk} \aaa'(T) = \lim_{\eps \lar 0} \int_0^T \int_0^t
f_{\eps}'(B_t - B_s) ds dt
\ee provided the limit exists. Formally, we can write

\be \aaa'(T) = \int_0^T \int_0^t \dd'(B_t - B_s) ds dt \ee Rosen was able to show that this integral
converges as $\eps \lar 0$, and proved an occupation time formula,
as well as some other facts about $\aaa'_t$. This paper deals with the derivative of ILT in 2 dimensions. In two dimensions,
we will use $f_{\eps}'$ to denote $\frac{\dd}{\dd x} f_{\eps}$. We let

\be \lll{nopk2} \aaa'_{\eps}(T) = \int_0^T \int_0^t
f_{\eps}'(B_t - B_s) ds dt \ee

Our main result is that $\aaa'_{\eps}(T)$ does not converge as $\eps \lar 0$. We will prove
that the asymptotic behavior as $\eps \lar 0$ is very similar to that which occurs for the ILT in 3 dimensions as
discovered by Yor. In particular, our main theorem is

\begin{theorem} \lll{bigone} $\Big{ \{}(\log(1/\eps))^{-1}
\aaa_{\eps}'(T), T \geq 0 \Big{ \}}$ converges in law to $\Big{ \{}
\frac{\sqrt{5}}{\pi 8 \sqrt[4]{2}}W_T, T \geq 0 \Big{ \}}$ as $\eps
\lar 0$, where $W_T$ is a one-dimensional Brownian motion.
\end{theorem}

Remark: $f_{\eps}'$ is an odd function, so $E[\aaa_{\eps}'(T)] = 0$,
which is why we need not subtract the expectation to obtain
convergence, as was required in the theorems of Yor, Varadhan, and Rosen.

We also will prove an analogous theorem about symmetric stable
processes. We let $X_t$ be a symmetric stable process of index $\bb$
with $1<\bb<2$, let $f_{t}$ be the density of $X_t$, and let $f'_t = \frac{\dd}{\dd x} f_t$. Again we will consider

\be \aaa_{\eps}'(T) = \int_0^T \int_0^t f_{\eps}'(X_{t}-X_{s})dsdt
\lll{a4}\ee and we will prove the following:

\begin{theorem} \lll{stab} $\{ \eps^{3/\bb - 3/2}\aaa_{\eps}'(T), T \geq 0 \}$ converges in law to $\{ c(\bb)W_T, T \geq 0 \}$ as $\eps \lar 0$ where $W_T$ is a one-dimensional
Brownian motion and $c(\bb)$ is given by

\be \frac{1}{2\sqrt{2}\pi^2} \Big( \int \int
\frac{1}{p^{\bb}}\frac{1}{q^{\bb}} \frac{1}{(p+q)^{\bb}}
e^{-(p^{\bb}+q^{\bb})}p_1 q_1 dpdq + \int \int \frac{1}{p^{2\bb}}
\frac{1}{(p+q)^{\bb}} e^{-(p^{\bb}+q^{\bb})}p_1 q_1 dpdq \Big) \ee
\end{theorem}

Included in the proof is the definition of the second integral in
the definition of $c(\bb)$. This integral does not converge
absolutely for $\bb \geq 3/2$, so we must clarify what it means.

The paper is organized as follows. Section 2 gives the outline of the proof. Sections 3 and 4 prove the bulk of the required
technical details. Section 5 wraps up the proof of Theorem \ref{bigone}. Section 6 deals with with the symmetric stable
case, and gives the proof of Theorem \ref{stab}.

\section{Outline of the proof}

The outline of the proof of Theorem \ref{bigone} follows closely the
proofs of Theorems 1, 2, and 3 given in \cite{rosen1}, though the
details, given in sections 3 and 4 below, are quite different. The
reader may refer to that paper to see a slightly different
presentation of the ideas of this section. We will show first that
the moments of $\aaa'_{\eps}(T)(\log(1/\eps))^{-1}$ converge to the
moments of a Brownian motion times $\frac{\sqrt{5}}{\pi
\sqrt{128\sqrt2}}$. Recall that

\be f_{\eps}(x) = \frac{1}{(2\pi \eps)^{d/2}} e^{\frac{-|x|^2}{2 \eps}} \ee We will express this in a form which is easier
for us to use:

\be \lll{psb} f_{\eps}(x) = \frac{1}{(2\pi)^{2}} \int e^{ipx-\eps p^{2}}
d^{2}p \lll{2}\ee This gives

\be \lll{psbn} f_{\eps}'(x) = \frac{i}{(2\pi)^{2}} \int p_1 e^{ipx-\eps p^{2}}
d^{2}p \lll{2}\ee We have

\be \lll{expect} E(\aaa'_{\eps}(T)^n) = \frac{i^n}{(2\pi)^{2n}}
\int_{({\bf R}^2)^n} \int_{D_T^n} e^{-\eps \sum_j p_j^2}
\prod_{j=1}^n p_{j,1} E\big[ \prod_{j=1}^n e^{i p_j(X_{t_j} -
X_{s_j})}\big] \prod_{j=1}^n ds_j dt_j d^2p_j \ee This is obtained
by combining $n$ copies of the integral (\ref{nopk2}) which defines $\aaa'_{\eps}(T)$, using the
definition (\ref{psb}) of $f_{\eps}$. Now, if $n$ is odd, then the
integrand is an odd function of $p$, and the expectation is
therefore 0. Since all odd moments of Brownian motion are $0$, we
need only show that the even moments converge to the right values.
The $2n$-th moment of Brownian motion at time $T$ is
$\frac{(2n)!}{2^n n!}T^n$, so we will show that

\be E[\aaa'_{\eps}(T)^{2n}] = \frac{(2n)!}{2^n
n!}\big(\frac{\sqrt{5}}{\pi
\sqrt{128\sqrt2}}(\log(1/\eps))\big)^{2n}T^n +
o(\log(1/\eps))^{2n}\ee This expectation is given by (\ref{expect})
with $n$ replaced by $2n$. The expectation on the right side of
(\ref{expect}) will depend on the ordering of the $s_j$'s and
$t_j$'s in $D_T^{2n}$. By independence, the expectation will factor
into the product of several expectations, each corresponding to a
component of the set $\unn_j [s_j,t_j]$. Each of these resulting
expectations are then integrated separately. Following
\cite{rosen1}, we will say that a component consisting of $m$
intervals $[s_j,t_j]$ is of order $m$. As in \cite{rosen1}, it turns
out that the dominant contribution comes from regions consisting of
$n$ components of order 2. Suppose, for the time being, that we are
integrating over only such regions. We will hold fixed the initial
points of the components $r_1 < ... < r_n$. In this case we will
show(Section 3) that each of the $n$ components contribute

\be \big(\frac{5}{\pi^2 128\sqrt2}(\log(1/\eps))\big)^{2} +
o(\log(1/\eps))^2 \ee The contribution of each configuration with
$n$ components of order 2 to (\ref{expect}) is therefore

\bea \lll{ned} \Big((\frac{5}{\pi^2 128\sqrt2}(\log(1/\eps))^{2} +
o(\log(1/\eps)^2)\Big)^n \int_{0 \leq r_1 \leq ... \leq r_n \leq T}
dr_1 ... dr_n \\ \nn = \frac{T^n}{2^n
n!}\big(\frac{\sqrt{5}\sqrt{2}}{\pi
\sqrt{128\sqrt2}}(\log(1/\eps))\big)^{2n} + o(\log(1/\eps))^{2n}
\\ \nn = \frac{T^n}{2^n
n!}\big(\frac{\sqrt{5}}{\pi 8 \sqrt[4]{2}}(\log(1/\eps))\big)^{2n} +
o(\log(1/\eps))^{2n} \eea where we have used the identity

\[ \int_{0 \leq r_1 \leq ... \leq r_n \leq T} dr_1 ...
dr_n=T^n/n!\] There are $(2n)!$ different ways to order the set
$\{s_1, ... , s_{2n}\}$, and each contributes (\ref{ned}), so the
total is

\be \frac{(2n)!}{2^n n!}\big(\frac{\sqrt{5}}{\pi
8\sqrt[4]{2}}(\log(1/\eps))\big)^{2n}T^n + o(\log(1/\eps))^{2n} \ee
which is, in fact, what we were aiming for. We must therefore show
that the contribution from all regions which have other than $n$
components of order 2 is $o(\log(1/\eps))^{2n}$. We will do this in
Section 4 by showing that each component of order $m$ with $m \geq
3$ contributes $o(\log(1/\eps))^m$. Note that any component of order
1, or indeed of any odd order, will in fact contribute 0, since the
integrand is an odd function. It will follow from all of this that
the moments of $(\log(1/\eps))^{-1} \aaa_{\eps}'(T)$ converge to the
right values. We will prove that the processes $(\log(1/\eps))^{-1}
\aaa_{\eps}'(T)$ are tight in Section 5, so that there is a limiting
process which has all of the same moments as $\frac{\sqrt{5}}{\pi
\sqrt{128\sqrt2}}W_t$. In addition, this process has independent
increments, as will also be shown in Section 5. These facts, taken
together, will prove Theorem \ref{bigone}.

\section{Components of order 2}

Here we will deal with the aforementioned components of order 2.
We'll begin by proving

\begin{proposition} \lll{order2} \be E[\aaa'_{\eps}(T)^{2}] = T \log (1/\eps)^2(\frac{10}{128\sqrt2
\pi^2} + o(1)) \ee \end{proposition}

To begin the proof, we write

\bea \lll{5} E[\aaa'_{\eps}(T)^{2}] = -\frac{1}{(2\pi)^{4}} \int
\int_{D_{T}^{2}} e^{-\eps(p^{2}+q^{2})} E[e^{ip(X_{t_{1}}-X_{s_{1}})
+ iq(X_{t_{2}}-X_{s_{2}})}]
\\ p_{1}q_{1}d^{2}sd^{2}tdpdq \nn \eea
where $D_{T}$ denotes the triangle $\{ (s_{i},t_{i})|0 \leq s_{i}
\leq t_{i} \leq T\}$. For simplicity we'll assume for the time being
that $T=1$. To handle the expectation in the integrand, we must
consider different orderings of the $s$'s and $t$'s. By symmetry, we
may assume $s_{1} < s_{2}$. We'll suppress the $\frac{-1}{(2\pi)^4}$
in front of the integral for the time being.

\medskip

{\bf Case 1:} $s_{1} < s_{2} < t_{2} < t_{1}$

\medskip

We rewrite the exponent in the expectation as
$ip(X_{t_{1}}-X_{t_{2}}) + i(p+q)(X_{t_{2}}-X_{s_{2}}) +
ip(X_{s_{2}}-X_{s_{1}})$, and then factor the expectation using
independence. As a result, the expectation becomes

\be e^{-p^{2}(a+c) - (p+q)^{2} b} \lll{6}\ee where $a =
t_{1}-t_{2}$, $b=t_{2}-s_{2}$, and $c=s_{2}-s_{1}$. Upon making this
linear transformation, the integral in question becomes

\be \int_{0}^{T} \Big[ \int \int \int_{a+b+c \leq t_{1}} dadbdc \int
\int e^{-p^{2}(a+c
+\eps)}e^{-(p+q)^{2}b}e^{-q^{2}\eps}p_{1}q_{1}dpdq \Big] dt_1
\lll{7}\ee

{\bf Case 2:} $s_{1} < s_{2} < t_{1} < t_{2}$

\medskip

We rewrite the exponent in the expectation as
$iq(X_{t_{2}}-X_{t_{1}}) + i(p+q)(X_{t_{1}}-X_{s_{2}}) +
ip(X_{s_{2}}-X_{s_{1}})$, and proceed as in Case 1. The integral in
question here is then

\be  \int_{0}^{T} \Big[ \int \int \int_{a+b+c \leq t_{1}} dadbdc
\int \int
e^{-p^{2}(c+\eps)}e^{-(p+q)^{2}b}e^{-q^{2}(a+\eps)}p_{1}q_{1}dpdq
\Big] dt_1 \lll{8}\ee where now $a = t_{2}-t_{1}$, $b=t_{1}-s_{2}$,
and $c=s_{2}-s_{1}$.

\medskip

{\bf Case 3:} $s_{1} < t_{1} < s_{2} < t_{2}$

\medskip

Here the expectation factors, and since the integrand of (\ref{5})
is then an odd function of $p$, the contribution to (\ref{5}) of
this case is 0.

\medskip

Remark: Before we begin, let us clear up a technical point that will
be necessary later. To do the computations, we will in fact
integrate da, db, and dc first. However, if we were to begin with dp
and dq, so that the integrals in cases 1 and 2 were of the form

\be \frac{-1}{(2\pi)^4}\int \int \int \int h(a,b,c,t) da db dc dt
\lll{9}\ee then h would be everywhere negative. Intuitively, this is
because when $p_1$ and $q_1$ are the same sign, $|p+q|$ is larger
than when they are opposite signs. In order to rigorously prove
this, just note that the map $\phi((p_1,p_2),((q_1,q_2)) =
((p_1,p_2),(-q_1,q_2))$ is a linear isometry which maps $U = \{p_1
q_1
>0\}$ bijectively onto $V = \{p_1 q_1 <0 \}$, and the function $|p+q|$ is greater at $(p,q) \in U$ than
at $\phi(p,q) \in V$. It will be pointed out later where we have
used this fact.

\medskip

We'll attack case 1 first. In order to compute the required
integral, we'll first examine the following:

\be  \int \int \int_{0 < a,b,c \leq 1} dadbdc \int \int
e^{-p^{2}(a+c+\eps)}e^{-(p+q)^{2}b}e^{-q^{2}\eps}p_{1}q_{1}dpdq
\lll{10}\ee Upon integrating $da$,$db$, and $dc$, we obtain

\bea \lll{11} \int \int \frac{(1-e^{-p^{2}})^{2}}{p^{4}}
\frac{(1-e^{-(p+q)^{2}})}{(p+q)^{2}} e^{-\eps(p^{2}+q^{2})}
p_{1}q_{1}dpdq \\ = \int \int \frac{(1-e^{-p^{2} /\eps})^{2}}{p^{4}}
\frac{(1-e^{-(p+q)^{2} /\eps})}{(p+q)^{2}} e^{-(p^{2}+q^{2})}
p_{1}q_{1}dpdq \nn \eea Upon converting to polar coordinates, with
$p=re^{i\th}$, $q=se^{i\phi}$, we arrive at

\bea \lll{12}\int \int \int \frac{(1-e^{-r^{2}
/\eps})^{2}}{r^{2}}s^{2} cos(\phi) e^{-(r^{2}+s^{2})} \\
\int_{0}^{2\pi} \frac{(1-e^{-|re^{i\th}+se^{i\phi}|^{2} /
\eps})}{|re^{i\th}+se^{i\phi}|^{2}}cos(\th) d\th d\phi drds \nn \eea
We'll isolate the $d\th$ integral. We may replace $\th$ with $\th +
\phi$. Note that

\be \frac{(1-e^{-|re^{i(\th+\phi)}+se^{i\phi}|^{2} /
\eps})}{|re^{i(\th+\phi)}+se^{i\phi}|^{2}} = \int_{0}^{1
/\eps}e^{-|re^{i(\th+\phi)}+se^{i\phi}|^{2}x} dx \lll{13}\ee so the
$d\th$ integral is

\bea \lll{14}\int_{0}^{1 / \eps} \int_{0}^{2\pi}
e^{-|re^{i(\th+\phi)}+se^{i\phi}|^{2}x} cos(\th+\phi) d\th dx \\
= \int_{0}^{1 / \eps} \int_{0}^{2\pi} e^{-|re^{i\th}+s|^{2}x}
[cos(\th)cos(\phi)-sin(\th)sin(\phi)] d\th dx \nn \eea We expand
this and pull the $cos(\phi)$ and $sin(\phi)$ terms outside of the
$d\th$ integral. Integrating the $sin(\phi)$ term with the
$cos(\phi)$ term already present outside the $d\th$ integral gives
0. Integrating the $cos(\phi)$ together with the other $cos(\phi)$
term already present gives $\pi /2$. We have therefore eliminated
$\phi$ from (\ref{12}). The contribution of (\ref{13}) to (\ref{12})
is therefore equal to

\be \frac{\pi}{2} \int_{0}^{1 / \eps} e^{-(r^{2}+s^{2})x}
\int_{0}^{2\pi} e^{-2rsx cos(\th)}cos(\th)d\th dx \lll{15}\ee By
\cite{gradrys} (p. 958, 8.431.5 with $\nu = 1$, replace $\th$ with
$\th + \pi$), what remains of the $d\th$ integral is equal to -$2
\pi I_{1}(2rsx)$, where $I_{1}$ denotes the modified Bessel function
of the first kind. Thus, (\ref{12}) is equal to

\bea \lll{16} -\pi^2 \int_{0}^{\infty} \int_{0}^{\infty} \int_{0}^{1
/ \eps}
\frac{(1-e^{-r^{2}/\eps})^{2}}{r^{2}}s^{2}e^{-(r^{2}+s^{2})(x+1)}I_{1}(2rsx)dxdrds
\\ = -\pi^2 \int_{0}^{\infty} \int_{0}^{\infty} \int_{0}^{ 1/ \eps}\int_{0}^{ 1/ \eps}
e^{-r^{2}y}(1-e^{-r^{2} /
\eps})s^{2}e^{-(r^{2}+s^{2})(x+1)}I_{1}(2rsx)dydxdrds \nn \eea
Expand this into two integrals, and rewrite the first one as

\be -\pi^2 \int_{0}^{\infty} \int_{0}^{ 1 / \eps}\int_{0}^{1 / \eps}
e^{-s^{2}(x+1)}s^{2}\int_{0}^{\infty}e^{-r^{2}(y+x+1)}I_{1}(2rsx)drdxdyds
\lll{17}\ee By \cite{gradrys} (p. 711,6.618.4) the $dr$ integral is
equal to

\be \frac{\sqrt{\pi}}{2
\sqrt{x+y+1}}e^{\frac{(2sx)^{2}}{8(x+y+1)}}I_{1/2}(\frac{(2sx)^{2}}{8(x+y+1)})
\lll{18}\ee This is $\frac{(e^{\frac{s^{2}x^{2}}{x+y+1}} -
1)}{\sqrt{2}sx}$, since $I_{1/2}(z)=\frac{1}{\sqrt{ 2\pi z}} (e^z -
e^{-z})$ by \cite{gradrys} (p.967 8.467). Thus, (\ref{16}) is

\bea \lll{19}\frac{-\pi^2}{\sqrt{2}}\int_{0}^{\infty} \int_{0}^{1 /
\eps}\int_{0}^{1 / \eps} e^{-s^{2}(x+1)}s
\frac{(e^{\frac{s^{2}x^{2}}{x+y+1}} - 1)}{x} dxdyds
\\= \frac{-\pi^2}{\sqrt{2}}\int_{0}^{1/\eps} \int_{0}^{1 /
\eps}\frac{1}{x}\int_{0}^{\infty} s
(e^{-s^{2}(x+1-\frac{x^{2}}{x+y+1})} - e^{-s^{2}(x+1)}) dsdydx \nn
\eea Since

\be \int_{0}^{\infty} se^{-bs^{2}} ds = (1/2)b^{-1} \lll{20}\ee we
see that (\ref{17}) is

\bea \lll{21}\frac{-\pi^2}{2\sqrt{2}} \int_{0}^{M}\int_{0}^{M}
\frac{1}{x}((x+1-\frac{x^{2}}{x+y+1})^{-1} - (x+1)^{-1})dxdy
\\=\frac{-\pi^2}{2\sqrt{2}} \int_{0}^{M}\int_{0}^{M}
\frac{x}{(x+1)(xy+2x+y+1)}dxdy \nn \eea where we have substituted
$M=1/\eps$ to simplify what follows. This last integral is
explicitly computable, but it is easier to calculate the derivative
and then apply L'Hospital's rule. To do so, we'll make use off the
following lemma:

\begin{lemma} \lll{fundy} If $h(x,y,M)$ is bounded, continuous in $x$ and $y$ on $\{ x,y \geq 0\}$, differentiable in $M$ with bounded derivative, then

\bea \lll{22}\frac{d}{dM} \int_{0}^{M} \int_{0}^{M} h(x,y,M) dxdy =
\int_{0}^{M} h(M,y,M)dy \\ + \int_{0}^{M} h(x,M,M) dx + \int_{0}^{M}
\int_{0}^{M} \frac{d}{dM}h(x,y,M) dxdy \nn \eea

\end{lemma}

In the case of (\ref{21}) the integrand does not depend on M, so the
last term is 0. We include the last term because it will be used
later. By the lemma, the derivative with respect to $M$ of
(\ref{21}) is

\bea \lll{23} \frac{-\pi^2}{2\sqrt{2}}\Big[ \frac{M}{M+1}
\int_{0}^{M} \frac{dy}{(M+1)y+2M+1} \\ \nn + \int_{0}^{M}
\frac{x}{(x+1)((M+2)x+M+1)} dx \Big]
\\ \nn =  \frac{-\pi^2}{2\sqrt{2}}\Big[ \frac{M}{M+1} \int_{0}^{M}
\frac{dy}{(M+1)y+2M+1} + \\ \nn \int_{0}^{M} \frac{1}{(x+1)} -
\frac{M+1}{((M+2)x+M+1)} dx \Big] \eea Performing the integration
gives

\bea \lll{hdj} \frac{-\pi^2}{2\sqrt{2}} \Big[ \frac{M}{(M+1)^{2}}
\log\Big( \frac{(M+1)M + 2M +1}{2M+1} \Big) + log(M+1)\\ \nn -
\frac{(M+1)}{(M+2)} \log \Big( \frac{(M+2)M + M+1}{M+1}\Big) \Big]
\eea The expression inside the brackets is asymptotic to
$\frac{2\log M}{M}$. This is an immediate consequence of the
following easily verified facts:

\medskip

1) $\frac{M}{(M+1)^{2}}= \frac{1}{M}+O(\frac{1}{M^{ 2}} )$

2) $\log\Big( \frac{(M+1)M + 2M +1}{2M+1} \Big) = \log M +O(1)$

3) $ log(M+1)= log(M)+O(\frac{1}{M} )$

4) $\log \Big( \frac{(M+2)M + M + 1}{M+1}\Big)= \log M
+O(\frac{1}{M} )$

5)
$\frac{(M+1)}{(M+2)}=1-\frac{1}{(M+2)}=1-\frac{1}{M}+O(\frac{1}{M^{
2}} )$.

\medskip

Thus, (\ref{23}) is equal to $\frac{\log
M}{M}(\frac{-\pi^2}{\sqrt{2}}+o(1))$, and it follows from
l'Hospital's rule that (\ref{21}) is equal to $(\log M)^2
(\frac{-\pi^2}{2\sqrt{2}}+o(1))$. Recall that we split (\ref{16})
into two integrals. We must now deal with the second, namely

\be -\pi^2 \int_{0}^{\infty} \int_{0}^{ 1 / \eps}\int_{0}^{1 / \eps}
e^{-s^{2}(x+1)}s^{2}\int_{0}^{\infty}e^{-r^{2}(y+x+1+1/\eps)}I_{1}(2rsx)drdxdyds
\lll{28} \ee We can follow steps (\ref{17})-(\ref{21}) exactly, with
the only difference being that we have $y+M$ in place of $y$. We get

\bea \lll{29} \frac{-\pi^2}{2\sqrt{2}} \int_{0}^{M}\int_{0}^{M}
\frac{1}{x}((x+1-\frac{x^{2}}{x+y+1 + M})^{-1} - (x+1)^{-1})dxdy
\\= \frac{-\pi^2}{2\sqrt{2}} \int_{0}^{M}\int_{0}^{M}
\frac{x}{(x+1)(xy+2x+y+1+M(x+1))}dxdy \nn \eea where, again,
$M=1/\eps$. We take the derivative as before, using Lemma
\ref{fundy}, and this time the integrand depends on $M$:

\bea \lll{30}\frac{-\pi^2}{2\sqrt{2}}\Big[\frac{M}{M+1} \int_{0}^{M}
\frac{dy}{(M+1)y+2M+1+M(M+1)}
\\+ \int_{0}^{M} \frac{x}{(x+1)((2M+2)x+2M+1)} dx \nn \\-
\int_{0}^{M}\int_{0}^{M} \frac{x}{(xy+2x+y+1+M(x+1))^2}dxdy \Big]
\nn \eea The first term

\be \frac{M}{M+1} \int_{0}^{M} \frac{dy}{(M+1)y+2M+1+M(M+1)}\ee

is bounded above by

\be \int_0^M \frac{dy}{M^2} = \frac{1}{M} \lll{31} \ee We may ignore
it, as it is $o(\log M /M)$. The second term

\bea \int_{0}^{M} \frac{x}{(x+1)((2M+2)x+2M+1)} dx \\ \nn =
\int_{0}^{M} \frac{1}{(x+1)} - \frac{2M+1}{((2M+2)x+2M+1)} dx\eea is
the same as the second integral in (\ref{23}), with $2M$ replacing
$M$. We can follow steps (\ref{23}) and (\ref{hdj}), and use fact
(3) above along with

6) $\log \Big( \frac{(2M+2)M + 2M + 1}{M+1}\Big)= \log M
+O(\frac{1}{M} )$

7)
$\frac{(M+1)}{(M+2)}=1-\frac{1}{(2M+2)}=1-\frac{1}{2M}+O(\frac{1}{M^{
2}} )$

\noindent to see that this term is asymptotic to $(1/2)\log (M)/M$.
The third term

\be \int_{0}^{M}\int_{0}^{M} \frac{x}{(xy+2x+y+1+M(x+1))^2}dxdy \ee
is bounded above by

\be \int_{0}^{M}\int_{0}^{M}\frac{dxdy}{M^2(x+1)^2} = O(1/M)
\lll{32} \ee and may also be ignored. Thus, (\ref{28}) is $(\log
(1/\eps))^2(\frac{-\pi^2}{8\sqrt{2}}+o(1))$. We have found the
asymptotics for (\ref{28}) and (\ref{17}). They are

\bea -\pi^2 \int_{0}^{\infty} \int_{0}^{ 1 / \eps}\int_{0}^{1 /
\eps}
e^{-s^{2}(x+1)}s^{2}\int_{0}^{\infty}e^{-r^{2}(y+x+1)}I_{1}(2rsx)drdxdyds
\lll{17q} \\ \nn = (\log (1/\eps))^2
(\frac{-\pi^2}{2\sqrt{2}}+o(1))\eea and

\bea -\pi^2 \! \! \int_{0}^{\infty} \! \! \int_{0}^{ 1 / \eps} \! \!
\int_{0}^{1 / \eps} \! \! e^{-s^{2}(x+1)}s^{2}\! \!
\int_{0}^{\infty} \! \! e^{-r^{2}(y+x+1+1/\eps)}I_{1}(2rsx)drdxdyds
\lll{28q} \\ \nn = (\log
(1/\eps))^2(\frac{-\pi^2}{8\sqrt{2}}+o(1))\eea Combining these as in
(\ref{16}), we see that (\ref{10}) is
$\log(1/\eps)^2(\frac{-3\pi^2}{8\sqrt{2}} + o(1))$

Now let us consider

\be  \int \int \int_{0 < a,b,c \leq t} dadbdc \int \int
e^{-p^{2}(a+c+\eps)}e^{-(p+q)^{2}b}e^{-q^{2}\eps}p_{1}q_{1}dpdq
\lll{34}\ee By a simple scaling we now show that this is equal to
$(\log (t/\eps))^2(\frac{-3\pi^2}{8\sqrt{2}}+o(1))$. The scaling is
as follows:

\bea \lll{35} \int \int \int_{0 < a,b,c \leq t} dadbdc \int \int
e^{-p^{2}(a+c+\eps)}e^{-(p+q)^{2}b}e^{-q^{2}\eps}p_{1}q_{1}dpdq
\\ \nn = \int \int \int_{0 < a,b,c \leq t} dadbdc \int \int
e^{-p^{2}t(\frac{a+c+\eps}{t})}e^{-(p+q)^{2}t\frac{b}{t}}e^{-q^{2}t(\frac{\eps}{t})}p_{1}q_{1}dpdq
\\ = t^2 \int \int \int_{0 < a,b,c \leq 1} dadbdc \int \int
e^{-p^{2}t(a+c+\eps/t)}e^{-(p+q)^{2}tb}e^{-q^{2}t(\eps/t)}(\sqrt{t}p_{1})(\sqrt{t}q_{1})dpdq
\nn \eea Now replace $(p,q)$ with $(p/\sqrt{t},q/\sqrt{t})$. The
$t^2$ in front of the integral is canceled, and we are left with

\bea \lll{36}\int \int \int_{0 < a,b,c \leq 1} dadbdc \int \int
e^{-p^{2}(a+c+\eps/t)}e^{-(p+q)^{2}b}e^{-q^{2}(\eps/t)}p_{1}q_{1}dpdq
\\ = (\log(t/\eps))^2(\frac{-3\pi^2}{8\sqrt{2}}+o(1))
\nn \eea We must now examine the same integral, but over the region
$\{a+b+c < t_1 \}$ rather than $\{0 < a,b,c \leq t_1\}$. However,
the remark following Case 3 shows that if $U \subs V$, then

\be \Big| \int_{a,b,c \in U} \Big| \leq \Big| \int_{a,b,c \in V}
\Big| \lll{37}\ee (Recall that both are negative) We also note that

\be \frac{(\log(t/\eps))^2}{(\log(1/\eps))^2} \lar 1 \lll{38} \ee We
can then write

\be \Big|\int_{a,b,c < (t_1 /3)}\Big| \leq \Big| \int_{a+b+c <
t_1}\Big| \leq \Big| \int_{a,b,c < t_1} \Big| \lll{39} \ee The first
and last integrals are both
$(log(1/\eps))^2(\frac{-3\pi^2}{8\sqrt{2}}+o(1))$, and it follows
that the middle one is as well. Thus, the integrand in the $dt_1$
integral (\ref{7}) is
$(\log(1/\eps))^2(\frac{-3\pi^2}{8\sqrt{2}}+o(1))$, with the $o(1)$
term uniformly bounded on $\{0< \dd < t_1 < T \}$. We will split up
 (\ref{7}) as:

\bea \lll{40} \int_{0}^{T} \! \int \! \int \! \int_{a+b+c \leq t_{1}} \!
dadbdcdt_{1} \! \int \! \int e^{-p^{2}(a+c
+\eps)}e^{-(p+q)^{2}b}e^{-q^{2}\eps}p_{1}q_{1}dpdq
\\= \int_0^\dd h(t_1,\eps) dt_1 + \int_\dd^T h(t_1,\eps) dt_1 \nn \eea
where $h$ denotes the result after doing the integrals in the other
variables. We know that, for the second integral,
$h(t_1,\eps)(\log(1/\eps))^{-2} \lar \frac{-3\pi^2}{8\sqrt{2}}$
uniformly. Thus, the second integral is
$(T-\dd)(\log(1/\eps))^2(\frac{-3\pi^2}{8\sqrt{2}}+o(1))$. The first
integral is bounded above in absolute value by(assuming $\dd<1$)

\bea \lll{41} \Big| \! \int_{0}^{\dd} \! \int \! \int \! \int_{a+b+c
\leq 1} \! dadbdcdt_{1} \! \int \! \int \! e^{-p^{2}(a+c
+\eps)}e^{-(p+q)^{2}b}e^{-q^{2}\eps}p_{1}q_{1}dpdq \Big|
\\= \dd (\log(1/\eps))^2(\frac{3\pi^2}{8\sqrt{2}}+o(1)) \nn \eea
By letting $\dd \lar 0$, we may finally conclude that the
contribution from Case 1 is

\be \lll{41a}T(\log(1/\eps))^2(\frac{-3\pi^2}{8\sqrt{2}}+o(1))\ee

Similar techniques will yield Case 2. We'll give only the outline
here. Recall that we are evaluating

\be \int_{0}^{T} \int \int \int_{a+b+c \leq t_{1}} dadbdcdt_{1} \int
\int
e^{-p^{2}(c+\eps)}e^{-(p+q)^{2}b}e^{-q^{2}(a+\eps)}p_{1}q_{1}dpdq
\lll{42} \ee As before, we begin by changing the domain to $\{0 \leq
a,b,c \leq 1\}$ and integrating $da$, $db$, and $dc$. We arrive at

\bea \lll{43}\int \int
\frac{(1-e^{-p^{2}})}{p^{2}}\frac{(1-e^{-q^{2}})}{q^{2}}
\frac{(1-e^{-(p+q)^{2}})}{(p+q)^{2}} e^{-\eps(p^{2}+q^{2})}
p_{1}q_{1}dpdq \\ = \int \int \frac{(1-e^{-p^{2}
/\eps})}{p^{2}}\frac{(1-e^{-q^{2} /\eps})}{q^{2}}
\frac{(1-e^{-(p+q)^{2} /\eps})}{(p+q)^{2}} e^{-(p^{2}+q^{2})}
p_{1}q_{1}dpdq \nn \eea We convert to polar coordinates again:

\bea \lll{44} \int \int \int (1-e^{-r^{2} /\eps})(1-e^{-s^{2}
/\eps}) cos(\phi) e^{-(r^{2}+s^{2})} \\ \int_{0}^{2\pi}
\frac{(1-e^{-|re^{i\th}+se^{i\phi}|^{2} /
\eps})}{|re^{i\th}+se^{i\phi}|^{2}}cos(\th) d\th d\phi drds \nn \eea
We follow the same steps for the $d\th$ integral as before(steps
(\ref{12})-(\ref{16})) to arrive at

\bea \lll{45} -\pi^2 \! \int_{0}^{\infty} \! \int_{0}^{\infty} \!
\int_{0}^{1 / \eps} \! (1-e^{-r^{2}/\eps})(1-e^{-s^{2}/\eps}) \\ \nn
e^{-(r^{2}+s^{2})(x+1)}I_{1}(2rsx)dxdrds \eea
We'll expand this into
4 integrals and do each separately. The first is

\bea \lll{46} -\pi^2 \int_{0}^{\infty} \int_{0}^{1 / \eps}
e^{-s^{2}(x+1)} \Big( \int_{0}^{\infty} e^{-r^{2}(x+1)}I_{1}(2rsx)dr
\Big) dxds \eea The $dr$ integral, as in step (\ref{18}), is equal
to $\frac{(e^{\frac{s^{2}x^{2}}{x+1}} - 1)}{\sqrt2 sx}$, and we
obtain:

\bea \lll{47} \frac{-\pi^2}{\sqrt 2}\int_{0}^{\infty} \int_{0}^{1 /
\eps} e^{-s^{2}(x+1)}\frac{(e^{\frac{s^{2}x^{2}}{x+1}} - 1)}{sx}
dxds \lll{47} \eea We use the identity

\be s \int_0^{\frac{x^2}{x+1}} e^{s^2t} dt =
\frac{(e^{\frac{s^{2}x^{2}}{x+1}} - 1)}{s} \lll{48} \ee And the
integral in question becomes

\bea \frac{-\pi^2}{\sqrt 2} \int_{0}^{1 / \eps}
\frac{1}{x}\int_{0}^{\frac{x^2}{x+1}} \int_{0}^{\infty}
e^{-s^{2}(x+1-t)}sdsdtdx \lll{49} \eea Note that $t \leq
\frac{x^2}{x+1} < x+1$, which implies $(x+1-t) > 0$, so there's no
problem with convergence. Tackling this integral again reduces to
basic calculus. Begin by substituting $u=s^{2}(x+1-t)$ to get

\bea \lll{50} \frac{-\pi^2}{2\sqrt 2}\int_{0}^{1 / \eps}
\frac{1}{x}\int_{0}^{\frac{x^2}{x+1}}
\frac{1}{x+1-t}\int_{0}^{\infty} e^{-u}dudtdx
\\ = \frac{-\pi^2}{2\sqrt 2}\int_{0}^{1 / \eps} \frac{1}{x}\int_{0}^{\frac{x^2}{x+1}}
\frac{1}{x+1-t}dtdx \nn
\\ = \frac{-\pi^2}{2\sqrt 2} \int_{0}^{1 / \eps} \frac{1}{x} log \Big(
\frac{x+1}{x+1-x^2 / (x+1)}\Big) dx \nn \eea If $M=1/\eps$, then, by
the fundamental theorem of calculus, $\frac{d}{dM}$ of the above is

\bea \frac{-\pi^2}{2\sqrt 2} \frac{\log \Big( \frac{M+1}{M+1-M^2 /
(M+1)}\Big)}{M} \lll{51} \eea
Since

\be \log \Big( \frac{M+1}{M+1-M^2 / (M+1)}\Big) = \log(M)+O(1)
\lll{52} \ee we see that (\ref{51}) is $\frac{\log
M}{M}(\frac{-\pi^2}{2\sqrt 2}+o(1))$, and thus our original integral
(\ref{46}) is $(\log M)^{2}(\frac{-\pi^2}{4\sqrt 2}+o(1))$. Recall
that (\ref{45}) was divided into four integrals. The remaining three
give a contribution of

\bea \lll{53} -\pi^2 \int_{0}^{\infty} \int_{0}^{\infty} \int_{0}^{1
/ \eps} (-e^{-r^2 /\eps}-e^{-s^2 /\eps}+e^{-r^2
/\eps}e^{-s^2 /\eps} ) \\
\nn e^{-(r^{2}+s^{2})(x+1)}I_{1}(2rsx)dxdrds \eea The integral
corresponding to the first and second term will be identical, and so
if we follow steps (\ref{46}) through (\ref{50}) we get

\bea  \lll{55} \frac{-\pi^2}{2\sqrt 2} \int_{0}^{1 / \eps}
\frac{1}{x} log
\Big( \frac{x+1+1/\eps}{x+1+1/\eps -x^2 / (x+1+1/\eps)}\Big) dx \\
\nn + 2 \frac{\pi^2}{2\sqrt 2} \int_{0}^{1 / \eps} \frac{1}{x} log
\Big( \frac{x+1+1/\eps}{x+1+1/\eps-x^2 / (x+1)}\Big) dx \\ \nn =
\frac{\pi^2}{2\sqrt 2} \int_{0}^{1 / \eps} \frac{1}{x} log
\Big(1-\frac{x^2}{(x+1+1/\eps)^2}\Big) dx \\
\nn - 2 \frac{\pi^2}{2\sqrt 2} \int_{0}^{1 / \eps} \frac{1}{x} log
\Big(1-\frac{x^2}{(x+1)(x+1+1/\eps)}\Big) dx \eea Each of these
integrals is $O(log(1/\eps))$. To see this, note that the absolute
values of the integrands are bounded above by $g(x) = \frac{-1}{x}
\log (1-x^2)$ on $\{0<x<\frac{1}{2}\}$. $g(x)$ is bounded on this
interval, so the integral over $\{0<x<\frac{1}{2}\}$ is bounded. For
$\{x
> \frac{1}{2}\}$ the integrand is bounded by $\frac{-1}{x} \log
(\frac{1}{2})$, since $1/\eps > x$, and it follows from this that
the integral is $O(\log(1/\eps)$. This proves that (\ref{44}) is
$\log(1/\eps)^2 (\frac{-\pi^2}{4\sqrt 2} + o(1))$. The remainder of
the proof that the integral we began with in case 2, namely
(\ref{8}), is $\log(1/\eps)^2 (\frac{-\pi^2}{4\sqrt 2} + o(1))$ is
identical to the steps (\ref{37}) through (\ref{41a}).

Combining our work in cases 1 and 2, and reinserting the constant
$\frac{-1}{(2\pi)^4}$ which was suppressed throughout, we see that

\be E[\aaa'_{\eps}(T)^{2}] = T \log (1/\eps)^2(\frac{5}{128\sqrt2
\pi^2} + o(1)) \ee We assumed at the outset that $s_1<s_2$, so this
must be multiplied by 2 to obtain the correct answer. \qed

\begin{corollary} The contribution to (\ref{expect}) of any
component of order 2 is $\log (1/\eps)^2(\frac{5}{128\sqrt2 \pi^2} +
o(1))$. That is, if $r_j$ is the left endpoint of a component of
order 2, and $r_{j+1}$ is the maximal right endpoint(see
(\ref{ned})), then the integral over this region is $\log
(1/\eps)^2(\frac{5}{128\sqrt2 \pi^2} + o(1))(r_{j+1}-r_j)$.
\end{corollary}

\medskip

{\bf Proof:} Suppose that the component of order 2 is composed of
$[s_k,t_k]$ and $[s_{k'}, t_{k'}]$. Then the integral in question is

\be \lll{99} \frac{-1}{(2\pi)^4} \int \! \! \int_{s_k, t_k, s_{k'},
t_{k'} \in [r_j, r_{j+1}]} e^{-\eps(p^{2}+q^{2})}
E[e^{ip(X_{t_{k}}-X_{s_{k}}) + iq(X_{t_{k'}}-X_{s_{k'}})}]
\\ p_{1}q_{1}d^{2}sd^{2}tdpdq
\ee

\medskip

We can write $X_{t_{k}}-X_{s_{k}}$ as $(X_{t_{k}}-X_{r_j})
-(X_{s_{k}}-X_{r_j})$, and likewise for $X_{t_{k'}}-X_{s_{k'}}$.
Since $X_{r_j+t} - X_{r_j}$ is itself a Brownian motion we see that
(\ref{99}) is equal to (\ref{5}) with $D_T$ replaced by
$D_{r_{j+1}-r_j}$. \qed

\section{Components of order $n \geq 3$}

We now turn to components of order $n$, where $n \geq 3$. We'll show

\begin{proposition} \lll{order3} The contribution to (\ref{expect}) of any component of order $n \geq 3$ is $o(\log (1/\eps)^n)$.
\end{proposition}

This entire section is devoted to the proof of this proposition.
Suppose that a component of order $n$ is can be formed by a specific
arrangement of intervals $[s_1,t_1], ... , [s_n,t_n]$ corresponding
to the variables $p_1, ... , p_n$. Let $\{ r_1, ... , r_{2n}\}$ be a
relabeling of the $s_i$'s and $t_i$'s so that $0 := r_0 \leq r_1
\leq r_2 \leq ... \leq r_{2n}$. We split up the expectation in the
integrand by independence and change coordinates as we did in the
order 2 case. The contribution of this arrangement of intervals is
then given by the integral

\be \lll{big} \int e^{-\eps \sum p_i^2} \prod_i (p_i)_1 \Big(
\int_{\sum c_j < T} \prod_j e^{-u_j^2 c_j} \prod_j dc_j \Big)
\prod_i dp_i \ee where $c_j = r_j - r_{j-1}$ and, as before,
$(p_i)_1$ denotes the x-coordinate of $p_i$. Each $u_j$ is a linear
combination of $p_i$'s, and is ordered in the natural way. That is,
$u_1 = p_1$, $u_2=p_1 + p_2$, $u_3 = p_2$ or $p_1$ or $p_1+p_2+p_3$,
etc. For each $j$, either $u_{j} - u_{j-1} = p_i$ or $u_{j} -
u_{j-1} = -p_i$ for some $i$. In the first case we'll refer to $u_j$
as increasing (abbreviated as $u_j \upa$) and in the second case
we'll say $u_j$ is decreasing (abbreviated as $u_j \doa$). To
simplify some of the notation that follows, let $u_0 = u_{2n} = 0$.
We will suppose first that our component contains no isolated
intervals. That is, there does not exist an interval $[s_i,t_i]$
such that $t_k, s_k \notin (s_i,t_i)$ for all $k$. In this case, we
may obtain a sufficient bound for (\ref{big}) by replacing the
integrand with the absolute value of the integrand. We may then also
replace the region $\sum c_j < T$ with $0<c_j<T$ for all $j$ and use
the simple fact that

\be \lll{bound} \int_0^T e^{-u_j^2 c_j} dc_j \leq
\frac{k}{(1+|u_j|)^2} \ee for a constant $k$ to reduce our problem
to bounding

\be \lll{abs} \int \frac{e^{-\eps \sum p_i^2} \prod
|p_i|}{\prod_{j=1}^{2n-1}(1+|u_j|)^2}\prod dp_i \ee We will
eventually need the following lemmas:

\begin{lemma} \lll{bnd1} \be \int_{\reals^2} \frac{e^{-\eps u^2}}{(1+|u|)^2} d^2 u =
O(\log(1/\eps))\ee

\end{lemma}

{\bf Proof:} Replace $u$ with $u/\sqrt{\eps}$ and the integral
becomes

\be \int_{\reals^2} \frac{e^{- u^2}}{(\sqrt{\eps} + u)^2} d^2 u \ee
The integral over $|u|>1$ is clearly $O(1)$, and for $|u|<1$ it
suffices to bound

\be \int_{|u|<1} \frac{1}{(\sqrt{\eps} + u)^2} d^2 u \ee Converting
to polar coordinates this is

\be 2\pi \int_0^1 \frac{r}{(\sqrt{\eps} + r)^2} dr \leq 2\pi
\int_0^1 \frac{1}{\sqrt{\eps} + r} dr = O(\log(1/\eps)) \ee

\qed

\begin{lemma} \lll{bnd2} There is a constant $c< \infty$ such that, independently
of $\eps, a$ and $n \geq 1$, we have

\be \int_{\reals^2} \frac{e^{- \eps u^2}}{(1 + |u|)^2
(1+|a+u|)^n}d^2 u < c \ee

\end{lemma}

{\bf Proof:} We can ignore the numerator. By H\"older's inequality
the integral is bounded by

\be \Big(\int \frac{1}{(1+|u|)^3} d^2 u\Big)^{2/3}\Big( \int
\frac{1}{(1+|u+a|)^{3n}} d^2 u\Big)^{1/3} < c < \infty \ee \qed

\bigskip

These lemmas motivate the intuition behind our approach, which we
first describe informally. The lemmas essentially say that a square
in the denominator gives a $\log$, whereas a cube or higher gives a
constant. We will write (\ref{abs}) in such a way that we can cancel
the $\prod |p_i|$ in the numerator with powers in the denominator.
We will use the Cauchy-Schwarz inequality to cut down on the number
of different terms in the numerator, and we will change variables by
a linear transformation. After all of this we will obtain a product
of a collection of integrals of the form in Lemma \ref{bnd1} with at
least one integral in the form of Lemma \ref{bnd2}. Each one of the
type in Lemma \ref{bnd1} contributes a $\log(1/\eps)$, whereas the
Lemma \ref{bnd2} type doesn't. When we multiply everything out, the
power of $\log(1/\eps)$ will be less than $n$. (As a side note,
Lemma \ref{bnd1} also indicates why we are initially not considering
the case of the isolated intervals. In that case there is some $p_i$
which is only present as a term in one $u_j$, so that if we were to
put the absolute value inside the integral as we are doing here, we
would have only a square of $p_i$ in the denominator with a $|p_i|$
present in the numerator. This would essentially give us

\be \int \frac{e^{-\eps p_i^2}}{1+|p_i|} d^2 p_i \ee And this is
only $O(\sqrt{1/\eps})$, which is not good enough.)

\medskip

To make good on this approach, we will need a way to make sure that,
after we cancel the terms in the numerator, we have enough terms
left in the denominator to obtain adequate convergence. In terms of
the sheer number of powers in the denominator there is no problem.
Lemma \ref{bnd1} suggests we need more than $2n$ powers on the
bottom, but there are $n$ powers on top versus $4n-2$ on the bottom,
for a total of $3n-2$ on the bottom. This is enough since $n \geq
3$. The tricky part is making sure that we have a proper assortment
on the bottom. The details are rather involved, so we first will
prove several technical lemmas. To state the first lemma, we need
another bit of terminology: We will say that $p_i$ is {\it t-free}
if there is no $t_k$ contained in $(s_i,t_i)$. Note that if $p_i$ is
t-free and $[s_i,t_i]$ is not an isolated interval, then at least
one $s_k$ is contained in $(s_i,t_i)$. We have the following lemma,
which was proved in \cite{rosen2}.

\begin{lemma} \lll{tless} The span of the decreasing $u_j$'s is equal to the span of the set of
all $p_i$'s which are not t-free. Furthermore, suppose that for each
t-free $p_i$ we chose $u(p_i)$ to be any one of the increasing
$u_j$'s which contains $p_i$ as a term. Then, if we let $D = \{$set
of decreasing $u_j$'s$\} \unn \{$set of all $u(p_i)$'s$\}$, $D$
spans the entire set $\{ p_1, ..., p_n \}$ \end{lemma}

{\bf Proof:} To begin with, suppose that $p_i$ is t-free. Then $p_i$
only appears as a term in increasing $u_j$'s, and so is not in the
span of the decreasing $u_j$'s. Conversely, suppose we have a
configuration of intervals such that the set of decreasing $u_j$'s
does not span the set of all $p_i$'s which are not t-free. Let
$p_{i_o}$ be the non-t-free $p$ with largest $s$ value which is not
in the span of the decreasing $u_j$'s. That is, if $s_i
>  s_{i_o}$ and $p_i$ is not t-free, then $p_i$ is in the span of
decreasing $u_j$'s. Now, let $u_{j_o}$ be the $u$ with largest $j$
value which contains $p_{i_o}$, i.e. such that $u_{j_o} - u_{j_o+1}
= p_{i_o}$. Then $u_{j_o+1}$ is decreasing, and we will obtain a
contradiction if we can show that $u_{j_o}$ is in the span of
decreasing $u_j$'s. If $u_{j_o}$ is decreasing there is nothing to
prove, so suppose $u_{j_o}$ is increasing.  Let $v>0$ be chosen as
small as possible so that $u_{j_o -v}$ is decreasing. The fact that
$p_{i_o}$ is not t-free implies that $p_{i_o}$ appears as a term in
$u_{j_o-v}$. Furthermore, we can write $p_{i_o} = u_{j_o-v} +
(u_{j_o-v+1} - u_{j_o-v}) + ... + (u_{j_o} - u_{j_o-1}) -
u_{j_o+1}$. Now, $u_{j_o-v}$ and $u_{j_o+1}$ are decreasing, and the
terms $(u_{j_o-v+1} - u_{j_o-v}), ... , (u_{j_o} - u_{j_o-1})$ are
all equal to $p_i$'s which have the properties that (i) they are not
t-free, because they appear as a term in $u_{j_o +1}$, and (ii) they
have larger $s$ values than $p_{i_o}$. We conclude that they are in
the span of the decreasing $u_j$'s, which means that $p_{i_o}$ is as
well. This is a contradiction, and establishes the first part of the
lemma. To prove the second part, just note that $u(p_i)$ contains
$p_i$ as a term as well as several other $p_k$'s which cannot be
t-free. These $p_k$'s are in the span of $D$ then, and thus $p_i$ is
as well. \qed

\medskip

In order to state the next lemma, we must consider (\ref{abs})
again. Let $u_{j_1}, ... , u_{j_n}$ be the increasing $u$'s in
order. That is, $u_{j_i} - u_{j_i - 1} = p_i$. We see that
(\ref{abs}) is bounded by

\be \int \lll{nabs} \frac{e^{-\eps \sum p_i^2} \prod (|u_{j_i}| +
|u_{j_i - 1}|)}{\prod_{j=1}^{2n-1}(1+|u_{j}|)^2}\prod dp_i \ee
Expand the numerator completely, and break this integral into the
sum of many integrals, each of which we do individually. Each of
these integrals has a product of $|u|$'s in the numerator, but no
$u$ can appear more than twice. This allows us to cancel all of the
$u$'s in the numerator with $u$'s in the denominator(Note: The word
"canceling", in this context, means replacing $\frac{|u|}{1+|u|}$
with 1). We arrive at the following integral:

\be \int \lll{abs2}\frac{e^{-\eps \sum
p_i^2}}{\prod_{j=1}^{2n-1}(1+|u_j|)^{m_j}}\prod dp_i \ee where $m_j
= 0,1,$ or $2$, depending on what power of $u_j$ appeared in the
numerator. The following lemma relates $m_j$ with the properties of
$u_j$ in the configuration of intervals.

\begin{lemma} \lll{facts} 1. If $u_j \doa$ then $m_j, m_{j-1} \geq 1$.

2. If $u_j \doa$ and $m_j=1$, then $u_{j+1} \upa$ and $m_{j+1} \geq
1$.

3. If $u_j$, $u_{j+1} \doa$ then $m_j = 2$.

\end{lemma}

{\bf Proof:} Each term in the numerator is of the form $(|u_j| +
|u_{j-1}|)$ where $u_j \upa$. We see that we can only have $m_j=0$
if $u_j$ appears in two terms in the numerator, and this can only
happen if both $u_j$ and $u_{j+1}$ are increasing. This proves (1).
If $u_j \doa$ then $u_j$ appears at most once in the numerator in
the term $(|u_{j+1}| + |u_j|)$ where $u_{j+1}$ must be increasing.
Furthermore, $u_{j+1}$ can appear in at most one other term, and so
if $m_j=1$ then $m_{j+1} \geq 1$. This proves (2). As for (3), if
$u_j$, $u_{j+1} \doa$ then $u_j$ does not appear in the numerator at
all, so $m_j = 2$.  \qed

\medskip

We now turn our attention to (\ref{abs}). Suppose that we can form
sets $A = \{a_1, ... , a_r \}$ and $B = \{b_1, ... , b_s \}$ with
the following properties:

\medskip

i) Each of $a_1, ... , a_r$ and $b_1, ... , b_s$ are equal to some
$u_j$.

ii) $A$ and $B$ each span $\{p_1,..., p_n\}$.

iii) If $a_i = u_j$ or $b_i = u_j$, then $m_j \geq 1$.

iv) If $a_i = b_k = u_j$, then $m_j=2$.

\medskip

Note that if we can find such sets we can, simply by deleting
elements if necessary, find two sets $A_o$ and $B_o$ which satisfy
the above properties and which each have $n$ elements. So in the
calculations which follow we'll assume that $n=s=r$, even though
when we eventually construct $A$ and $B$ they may have more than $n$
elements. Given two such sets, we could bound (\ref{abs2}) by

\be \int \lll{abs3}\frac{e^{-\eps \sum
p_i^2}}{K(a_1,...,a_n)\prod_{j=1}^{n}(1+|a_j|)(1+|b_j|)}\prod dp_i
\ee where $K$ is of the form $(1+|u_j|)$ for some $j$. Recall that
we have $3n-2$ powers of $u$'s in the denominator, so there will
always be at least one term left over after choosing our sets $A$
and $B$. This term will contain a linear combination of $p_i$'s, but
since $A$ spans $\{p_1,..., p_n\}$ we may write it as a linear
combination of $a_j$'s. It is irrelevant what the linear combination
present in $K$ actually is, except that it must be nonzero, of
course. Now, we can apply the Cauchy-Schwarz inequality to bound
(\ref{abs3}) by

\bea \Big( \int \lll{abs4} \frac{e^{-\eps \sum
p_i^2}}{K(a_1,...,a_n)^2 \prod_{j=1}^{n}(1+|a_j|)^2}\prod dp_i \Big)
^{1/2} \\ \times \nn \Big( \int \frac{e^{-\eps \sum p_i^2}}{
\prod_{j=1}^{n}(1+|b_j|)^2}\prod dp_i \Big) ^{1/2}\eea There is a
constant $c>0$ so that $\sum p_i^2 > c \sum a_i^2$ and $\sum p_i^2 >
c \sum b_i^2$; this is because the functions $\sum a_i^2$ and $\sum
b_i^2$ are homogeneous of degree 2 in the $p_i$'s and bounded on
$\sum p_i^2 = 1$. Thus, (\ref{abs4}) is bounded by

\bea \Big( \int \lll{abs5} \frac{e^{-\eps c \sum
a_j^2}}{K(a_1,...,a_n)^2 \prod_{j=1}^{n}(1+|a_j|)^2}\prod dp_i \Big)
^{1/2} \\ \times \nn \Big( \int \frac{e^{-\eps c \sum b_j^2}}{
\prod_{j=1}^{n}(1+|b_j|)^2}\prod dp_i \Big) ^{1/2}\eea Now, we apply
a linear change of coordinates to these integrals so that we are
integrating with respect to $a_j$ and $b_j$ instead of $p_i$.
Relabel if necessary so that $a_1$ is one of the $a$'s which appears
as a term in $K(a_1,...,a_n)$. We see that the first integral in
(\ref{abs5}) is bounded by a constant times

\be \lll{abs6} \int \Big(\int \frac{e^{-c \eps
a_1^2}}{K(a_1,...,a_n)^2(1+|a_1|)^2}\,da_1 \Big)
\prod_{j=2}^{n}\frac{e^{-c \eps a_j^2}}{(1+|a_j|)^2}da_j \ee By
Lemma \ref{bnd2} the inner integral is $O(1)$ and by Lemma
\ref{bnd1} the others are all $O(\log(1/\eps))$. Lemma \ref{bnd1}
also shows that the second integral in (\ref{abs5}) is $O(\log^n
(1/\eps))$. We see that (\ref{abs5}) is
$O((\log(1/\eps))^{n-(1/2)}$, and this shows that (\ref{big}) is
$o((\log(1/\eps))^n)$, which is what we set out to prove.

All that remains, then, is to show that we can always find sets $A$
and $B$ of this form. For this, we'll use Lemmas \ref{tless} and
\ref{facts}. Lemma \ref{tless} gives us a good first initial
candidate for $A$ and $B$. We can let $A$ be equal to the set of
(distinct) decreasing $u_j$'s together with elements $u(p_i)$ for
each t-free $p_i$ (Recall that all decreasing $u_j$'s have $m_j \geq
1$). The possible problem with this is that every increasing $u_i$,
and in particular each possibility for $u(p_i)$, appears at least
once in the numerator of (\ref{nabs}), so that we need to make sure
that we really can appropriately choose the $u(p_i)$'s.
Nevertheless, as will be shown below, this works for $A$. $B$ cannot
be chosen the same way, however. This is because if $u_j$ is
decreasing but $u_{j+1}$ is increasing, then $u_j$ appears exactly
once in the numerator of (\ref{nabs}) and we may have $m_j = 1$, so
that $u_j$ cannot be in both $A$ and $B$. $B$ will have to be formed
in a different manner.

In order to eventually create the set $B$, we'll begin by creating
an increasing collection of sets $B_n$ by considering decreasing
$u_j$'s for increasing values of $j$. Start with the smallest $j$
such that $u_j$ is decreasing. If $m_j = 2$, then let $B_1=\{u_j\}$.
If $m_j = 1$, then we know from Lemma \ref{facts} that $u_{j+1}$ is
increasing. We can then choose $d > 0$ such that $u_j = u_{j+1} -
(u_{j+d} - u_{j+d+1})$; $d$ is simply chosen to be the largest value
such that $u_{j+d}$ contains $p_i = u_{j+1} - u_j$ as a term. let
$B_1 = \{u_{j+1}, u_{j+d}, u_{j+d+1}\}$. We will essentially repeat
this for each decreasing $u_j$. Suppose that the set $B_n$ has
already been formed. Let $j$ be as small as possible so that $u_j
\notin span\{ B_n \}$ and $u_j$ is decreasing. If $m_j = 2$, then
let $B_{n+1}=B_n \unn \{u_j\}$. If $m_j = 1$, then let $B_{n+1} =
 B_n \unn \{u_{j+1}, u_{j+d}, u_{j+d+1}\} \setminus \{u_j \}$,
where again $d$ is such that $u_{j+d+1} - u_{j+d} = u_{j+1} - u_j$.
The reason for subtracting the set $\{ u_j \}$ is that it may
already be in the set $B_n$, having having been of the form
$u_{j'+d'}$ or $u_{j'+d'+1}$ for an earlier $j'$. Repeat this
process through all of the decreasing $u_j$'s. The final set
obtained, say $B_N$, will span the set of decreasing $u_j$'s. To see
this, suppose to the contrary, and let $u_{j_o}$ be the $u_j$ with
largest $j$ value which is not in the span of $B_N$. Clearly then
$m_{j_o}=1$, which means that an earlier $B_n$ must have contained
$u_{j_o+1}, u_{j_o+d_o}$, and $u_{j_o+d_o+1}$. Any of these elements
which are increasing must be present in $B_N$, and any decreasing
ones have larger $j$ values than $u_{j_o}$, which means they are in
the span of $B_N$. Thus, $u_{j_o} = u_{j_o+1} - (u_{j_o+d_o} -
u_{j_o+d_o+1})$ is also in the span of $B_N$, a contradiction. The
set $B_N$ also satisfies property (iii) above. This will be shown
using the following lemma:

\begin{lemma} \lll{bfacts}1. If $u_i$ is in $B_N$ then $u_i$ is either decreasing or else
neighbors on a decreasing interval(i.e. at least one of $u_{i-1}$
and $u_{i+1}$ is decreasing).

2. If $u_i, u_{i+1}$ are both increasing and $u_i \in B_N$ then
$u_{i-1} \doa$, $m_{i-1} = 1$, and $m_i \geq 1$.

\end{lemma}

{\bf Proof:} If $u_i$ is increasing and in $BN$ then $u_i$ must be
of the form $u_{j+1}, u_{j+d},$ or $u_{j+d+1}$ for some $j$ where
$j$ and $d$ are as in the construction of the set $B'$ above, i.e.
$u_j$ is decreasing, $m_j = 1$, and $u_{j} = u_{j+1} - (u_{j+d} -
u_{j+d+1})$. It is always true that $u_{j+d+1}$ is decreasing, so
this cannot be $u_i$. (1) is proved by noting that $u_{j+1},u_{j+d}$
are neighbors to the decreasing intervals $u_j, u_{j+d+1}$
respectively. If, in addition, the situation in (2) arises then
$u_i$ cannot be of the form $u_{j+d}$ since in that case $u_{j+d+1}$
would be decreasing. Thus, $u_i$ is of the form $u_{j+1}$. In order
for $u_{j+1}$ to be included in $B'$ it was necessary that $u_j
\doa$ and $m_j=1$. By part 2 in Lemma \ref{facts} $m_i \geq 1$. \qed

If $u_j \in B_N$ then either $u_j \doa$, in which case $m_j \geq 1$
by part 1 of Lemma \ref{facts}, or $u_j \upa$ in which $m_j \geq 1$
by part 1 of Lemma \ref{facts} or part 2 of Lemma \ref{bfacts},
depending on whether $u_{j+1} \doa$ or $\upa$. Thus, $B_N$ satisfies
(iii) as claimed. Let $B' = B_N$ and $A'$ be the set of all
decreasing $u_j$'s. We know from the discussion above that $B'$
satisfies (i) and (iii). $A'$ clearly satisfies (i), and satisfies
(iii) by part 1 of Lemma \ref{facts}. $A'$ and $B'$ together satisfy
(iv) because of the way that $B'$ was constructed, and both span the
set of all decreasing $u_j$'s. We need now only extend them to sets
$A$ and $B$ which span all of $\{p_1,..., p_n\}$. $A'$ and $B'$
already span the set of all non-t-free $p_i$'s, by the first part of
Lemma \ref{tless}. In light of the second part of Lemma \ref{tless},
all that remains is to show that, for any t-free $p_i$, we can
always choose $u_1(p_i), u_2(p_i)$ which contain $p_i$ as a term,
and which we may include in $A$ and $B$ respectively without
violating rules (iii) and (iv).

Suppose $p_i$ is t-free, and $k$ is chosen as large as possible so
that $s_i < s_{i+1} < ... < s_{i+k} < t_i$(Note: $k$ here is not the
same as $d$ above; they differ by 1). Let $u_{j} - u_{j-1} = p_i$.
The term $p_i p_{i+1}...p_{i+k}$ in (\ref{big}) becomes
$(|u_{j}|+|u_{j-1}|)...(|u_{j+k}|+|u_{j+k-1}|)$ in (\ref{nabs}),
with $u_j, ... , u_{j+k}$ not appearing anywhere else in the
numerator. If $k>1$ we can just note that, upon expanding this
expression, the sum of the powers of $u_{j+k}$ and $u_{j+k-1}$ in
the numerator is at most two. This means that $m_{j+k} + m_{j+k-1}$
must be at least $2$, and we can choose $u_1(p_i), u_2(p_i)$ as some
combination of $u_{j+k}$ and $u_{j+k-1}$. It is possible that
$u_{j+k}$ is already in $B'$, and so we must make sure that if
$u_1(p_i) = u_{j+k} \neq u_2(p_i)$ that we interchange $u_1(p_i)$
and $u_2(p_i)$, so that $u_{j+k}$ is not in both $A$ and $B$, so as
to not violate (iv). Note that if $k>1$ then $u_{j+k-1} \notin B'$
by Lemma \ref{bfacts}, since $u_{j+k-1} \upa$ and neighbors only on
increasing intervals. In the case that $k=1$ we still have $m_j +
m_{j+1} \geq 2$, but now it is possible that both $u_j$ and
$u_{j+1}$ are in $B'$, since both neighbor upon intervals which may
be decreasing. However, if this is the case then, since $u_j,
u_{j+1} \upa$, we have by Lemma \ref{bfacts} $u_{j-1} \doa$ and
$m_{j-1}=1$. Recall that we have the term
$(|u_{j}|+|u_{j-1}|)(|u_{j+1}|+|u_j|)$ in the numerator, with
$u_{j-1}, u_j , u_{j+1}$ appearing nowhere else in the numerator.
The sum of the powers of $u_{j-1}, u_{j},$ and $u_{j+1}$ is two, and
thus $m_{j-1}+m_j + m_{j+1} = 4$. Since $m_{j-1} = 1$, one of $m_j$
and $m_{j+1}$ is 2. We can then let $u_1(p_i) = u_2(p_i) = u_j$ or
$u_{j+1}$, depending on whether $m_j$ or $m_{j+1}$ is 2. This
handles the case $k=1$. (If $k=0$ then we would have an isolated
interval, and this argument doesn't work. This is the only place
where we used the fact that we had no isolated intervals.) Doing
this for each t-free $p_i$ we create the sets $A$ and $B$, which are
guaranteed by Lemma \ref{tless} to satisfy the property (ii). $A$
and $B$ also satisfy properties (i), (ii), and (iv) by construction,
so we have completed the proof in the case where no isolated
intervals are present.

Now for the isolated intervals case. Recall that the integral which
gives us the contribution from this configuration is

\be \lll{big2}  \int e^{-\eps \sum p_i^2} \prod
(p_i)_1\prod_{j=1}^{2n-1}\Big( \int_{\sum t_j < T} \prod_j e^{-u_j^2
t_j} \prod_j dt_j\Big) \prod dp_i \ee As mentioned before, here we
cannot replace the integrand with its absolute value, for in that
case each isolated interval would contribute a $\sqrt{1/\eps}$ to
the integral. Cancelation occurs in the integral, however, since the
integrand is positive in some regions and negative in others. It
turns out that it is enough to integrate each of the variables
corresponding to isolated intervals first, and then to bring the
absolute value inside the integral. After we have "removed" the
initial set of isolated intervals in this fashion, we will have
created a new configuration of intervals, which may again contain
isolated intervals. We can remove these isolated intervals by a
different method than was used for the first set. This brings us to
a new configuration, which may again have isolated intervals, which
we again remove, etc. After a finite number of steps we either have
removed all intervals or we have arrived at an arrangement with no
isolated intervals. In the second case we are reduced to the case we
have already done, and the first is handled easily in a slightly
different way.

Let us bring in some definitions in order to make this rigorous. Let
our initial configuration of intervals be denoted $K_0$, and let
$K_m$ be the configuration of intervals obtained upon removing the
isolated intervals from $K_{m-1}$. We will say $p_i \in K_m$ to mean
that the interval $(s_i,t_i)$ appears in the configuration $K_m$,
and we will define the {\it order} of $K_m$ to be the number of
$p_i$'s in $K_m$. Let $u_{(m,1)}, ... , u_{(m,n_m)}$ be the linear
combinations of $p_i$'s which appear in the configuration $K_m$,
ordered from left to right. Let us define $I_m$ to be the set of all
$j$ values corresponding to isolated intervals in $K_m$; that is,
$I_m=\{j : u_{(m,j)} = p_i$ where $(s_i,t_i)$ contains no $s_k$ or
$t_k$ in $K_m \}$. A $\hat{p}$ will refer to the $p$ associated to
an isolated interval. That is, if $j \in I_m$ and $p_i$ is the $p$
which appears only in $u_{(m,j)}$, label $p_i$ as $\hat{p}_{m,j}$.
We can bound (\ref{big2}) by

\be \lll{big3} \int \prod e^{-\eps p_i^2}|p_i| \Big( \int \prod_{j
\notin I_0} e^{-u_{(0,j)}^2 t_j} \Big|\prod_{j \in I_0} \! \! \int
\! \! \int \!
 e^{-\eps \hat{p}_{(0,j)}^2}  (\hat{p}_{(0,j)})_1 e^{-u_{(0,j)}^2
t_j} dt_j d\hat{p}_{(0,j)}\Big| \prod_{j \notin I_0} dt_j \Big)
\prod dp_i\ee The first and last products over all $i$ such that
$p_i \neq \hat{p}_{(0,j)}$ for all $j$. We will get a good bound on
the $dt_j d\hat{p}_j$ integrals. Note that we have suppressed the
region of integration in $t_j$, since it may be quite complicated.
We do know that the upper limit of integration is bounded above by
$T$, and this allows us to get a sufficient bound, as the following
lemma shows.

\begin{lemma} \lll{bnd3} For any $a$ with $0<a<T$ and any $k \in {\bf R}^2$, we have

\be \Big|\int \int_0^a e^{-\eps p^2} p_1 e^{-(p+k)^2 t} dt dp\Big| =
|k|O(log(1/\eps)) \ee independently of $a$.

\end{lemma}

{\bf Proof:}

\bea \int \int_0^a e^{-\eps p^2} p_1 e^{-(p+k)^2 t} dt dp \\
\nn = \int_0^a \Big( \int e^{-\eps p_1^2} p_1 e^{-(p_1+k_1)^2 t}
dp_1 \int e^{-\eps p_2^2} e^{-(p_2+k_2)^2 t} dp_2\Big) dt\eea Now,
for $p,k \in R^2$, $(p+k)^2 = p^2 + k^2 + 2 p \cdot k$, so this is

\bea \lll{hiya} \int_0^a e^{-k^2 t}  \Big( \int e^{-\eps p_1^2} p_1
e^{-(p_1^2+2 p_1 k_1) t} dp_1 \int e^{-\eps p_1^2} e^{(-p_2^2+2 p_2
k_2) t} dp_2\Big) dt \\
\nn = \int_0^a e^{-k^2 t} e^{\frac{k^2 t^2}{\eps + t}} \Big( \int
e^{-(\eps+t) (p_1+\frac{k_1 t}{\eps+t})^2} p_1 dp_1 \int
e^{-(\eps+t) (p_2+\frac{k_2 t}{\eps+t})^2} dp_2\Big) dt \\ \nn =
\int_0^a e^{-k^2 t} e^{\frac{k^2 t^2}{\eps + t}} \Big( \int
e^{-(\eps+t) p_1^2} (p_1- \frac{k_1 t}{\eps+t})dp_1 \int
e^{-(\eps+t) p_2^2} dp_2\Big) dt \eea We now split the $p_1$
integral into two pieces, and we see that the first one,

\be \int e^{-(\eps+t) p_1^2} p_1 dp_1 \ee is 0 by symmetry (this is
what will give us the extra convergence). We use the fact that, for
$d=1,2$ we have

\be \int e^{-(\eps+t) p_d^2} dp_d = \frac{c}{\sqrt{\eps+t}} \ee for
some constant $c$. We'll also replace $\frac{t}{\eps + t}$ and
$e^{-k^2 t} e^{\frac{k^2 t^2}{\eps + t}}$ by the trivial bound of 1.
This shows us that we can bound (\ref{hiya}) by

\be c^2 |k_1| \int_0^a \frac{1}{\eps + t}dt \ee Since $a<T$, this is
$|k|O(\log(1/\eps))$, independently of $a$. \qed

\medskip

We integrate the $\hat{p}_{0,j}$'s first, and by the previous lemma
each one gives $|u_{0,j-1}|O(\log(1/\eps))$ ($|u_{0,j-1}|$ is the
$u_j$ which appears immediately before and after the isolated
interval corresponding to $\hat{p}_{0,j}$). (\ref{big3}) is thus

\bea \lll{big4}  O(\log(1/\eps))^{|I_0|} \! \int \! \prod_{p_i \neq
\hat{p}_{(0,j)} \forall j} e^{-\eps p_i^2} |p_i| \prod_{(j+1) \in
I_0} |u_{(0,j)}|
\\ \nn \Big( \int \prod_{j \notin I_0} e^{-u_{(0,j)}^2 t_j} \prod_{j
\notin I_0} dt_j \Big) \prod_{p_i \neq \hat{p}_{(0,j)} \forall j}
dp_i \eea Since the integrand is now positive we can extend the
region of integration for the $t_i$'s to be $0 < t_i < T$ and use
(\ref{bound}) to bound (\ref{big4}) by

\bea \lll{big5}  O(\log(1/\eps))^{|I_0|} \! \int \! \! \prod_{p_i
\neq \hat{p}_{(0,j)} \forall j} e^{-\eps p_i^2}  |p_i| \prod_{(j+1)
\in I_0} \! \! \\ \nn |u_{(0,j)}| \prod_{j \notin I_0} \! \!
\frac{1}{1+u_{(0,j)}^2} \prod_{p_i \neq \hat{p}_{(0,j)} \forall j}
\! \! dp_i \eea Suppose that $u_{m,j}$ is an isolated interval in
$K_m$. Then $u_{m,j-1} = u_{m,j+1}$. We will say in this case that
$u_{m,j-1}$ contains $u_{m,j}$. If $u_{m,j} = u_{m',j'}$, where
$m>m'$, and $u_{m',j'+1}$ is isolated in $K_{m'}$, we will also say
that $u_{m,j}$ contains $u_{m',j'+1}$. We will let $l_{m,j}$ denote
the total number of isolated intervals which the interval $u_{m,j}$
contained in all $K_{m'}$'s, where $m'<m$. Each $u_{1,j}$ which
contained one or more isolated intervals in $K_0$ will appear to a
power $l_{1,j}$ in the numerator of (\ref{big5}) as a result of
Lemma \ref{bnd3}, but the term $(1+u_{1,j}^2)$ will also appear an
extra $l_{1,j}$ times in the denominator. We see that (\ref{big5})
is

\be \lll{big5a}  O(\log(1/\eps))^{|I_0|} \int e^{-\eps \sum_{p_i \in
K_1} p_i^2} \prod_{p_i \in K_1} |p_i| \prod_{1 \leq j \leq
n_1}|u_{(1,j)}|^{l_{1,j}} \frac{1}{(1+u_{(1,j)}^{2+2l_{1,j}})}
\prod_{p_i \in K_1} dp_i \ee We must have some idea how the integral
(\ref{big5}) can be bounded as we remove successive stages of
isolated intervals, and Lemma \ref{remove} below gives us that. The
following lemma prepares us to prove Lemma \ref{remove}.

\begin{lemma} \lll{bnd4}

\be \lll{hi} \int e^{-\eps p^2}|p| \frac{1}{(1+|k+p|)^m} dp =
(1+|k|)O(1) + O(\log(1/\eps)) \ee if $m=3$, and is $(1+|k|)O(1)$ if
$m>3$.

\end{lemma}

{\bf Proof:} (\ref{hi}) is bounded by

\bea \int e^{-\eps (p-k)^2}(|p|+|k|) \frac{1}{(1+|p|)^m} dp \eea
Divide this into two integrals. The one with $|k|$ in the numerator
is bounded by

\be |k| \int \frac{1}{(1+|p|)^{m}} dp = |k|O(1) \ee The other is
bounded by

\be c \int e^{-\eps (p-k)^2} \frac{1}{1+|p|^{m-1}}dp \ee Again if
$m>3$ this is $O(1)$. If $m=3$, divide the region into
$\{|p|>2|k|\}$ and $\{|p|<2|k|\}$. On $\{|p|>2|k|\}$ we can bound
the integral by

\be \int e^{\eps p^2/2} \frac{1}{1+p^2} = O(log(1/\eps))\ee by Lemma
\ref{bnd1}. On $\{|p|<2|k|\}$ we can bound it by

\be \int_{|p|<2|k|} \frac{1}{1+p^2} dp \leq \log(|k|+1) \leq |k| \ee
These bounds combine to prove the lemma. \qed

\begin{lemma} \lll{remove} Suppose that $K_m$ contains isolated intervals. Then (\ref{big2}) is

\bea \lll{big8}  O(\log(1/\eps))^{|I_0| + ... + |I_{m}|} \int
e^{-\eps \sum_{p_i \in K_{m+1}} p_i^2} \prod_{p_i \in K_{m+1}} |p_i|
\\ \nn\prod_{1 \leq j \leq n_{m+1}}(1+|u_{(m+1,j)}|)^{l_{m+1,j}}
\frac{1}{(1+u_{(m+1,j)}^{2+2l_{m+1,j}})} \prod_{p_i \in K_{m+1}}
dp_i \eea

\end{lemma}

\medskip

{\bf Proof:} By induction. We know that it is true for
$m=0$(see(\ref{big5a})). Assume that it is true for $m-1$, so
(\ref{big2}) is

\bea \lll{big6}  O(\log(1/\eps))^{|I_0| + ... + |I_{m-1}|} \int
e^{-\eps \sum_{p_i \in K_m} p_i^2} \prod_{p_i \in K_m} |p_i|
\\ \nn \prod_{1 \leq j \leq n_m}(1+|u_{(m,j)}|)^{l_{m,j}}
\frac{1}{(1+u_{(m,j)}^{2+2l_{m,j}})} \prod_{p_i \in K_m} dp_i \eea
We will integrate the variables in $K_m$ corresponding to isolated
intervals. We can rewrite the integral in (\ref{big6}) as

\bea \lll{big7} \int e^{-\eps \sum_{p_i \in K_m, i \notin \hat{I}_m}
p_i^2} \prod_{p_i \in K_m, i\notin \hat{I}_m} |p_i| \prod_{j \notin
I_m}\frac{(1+|u_{(m,j)}|)^{l_{m,j}}}{(1+u_{(m,j)}^{2+2l_{m,j}})}
\\ \nn \Big(\prod_{j \in I_m} \int |\hat{p}_{m,j}|e^{-\eps \hat{p}_{m,j}}
\frac{(1+|u_{(m,j)}|)^{l_{m,j}}} {(1+u_{(m,j)}^{2+2l_{m,j}})}
d\hat{p}_{m,j}\Big) \prod_{p_i \in K_m i\notin \hat{I}_m} dp_i \eea
It is simply to verify that

\be \frac{(1+|u_{(m,j)}|)^{l_{m,j}}} {(1+u_{(m,j)}^{2+2l_{m,j}})}
\leq K \frac{1} {(1+u_{(m,j)}^{2+l_{m,j}})}\ee For some constant $K$
depending on $l_{m+1,j}$. Each $\hat{p}_{m,j}$ integral is

\be (1+|u_{m,j} - \hat{p}_{m,j}|)O(\log(1/\eps)) \ee by Lemma
\ref{bnd4}. Plugging this into (\ref{big7}) and relabeling the $u$'s
with index $m+1$ instead of $m$ gives (\ref{big8}). \qed

To complete the proof of Proposition \ref{order3}, let us consider
several cases. Recall that {\it order} refers to how many intervals
$[s_i, t_i]$ make up a configuration.

\medskip

{\bf Case 1:} There is a $K_m$ of order greater than or equal to 3
which contains no isolated intervals.

\medskip

In this case our integral in (\ref{big8}) is almost the same as what
would have been obtained if we had started with the configuration
$K_m$. The only difference is the presence of the extra powers
$l_{m,j}$, which in fact cause greater convergence. Thus, by what we
did earlier in this section, the remaining integral is
$o(\log(1/\eps))^{|K_m|}$. Since

\be |I_0| + ... + |I_{m-1}| + |K_m| = n \ee we see that (\ref{big2})
is $o(\log(1/\eps))^n$, which is what we set out to prove.

\medskip

{\bf Case 2:} There is a $K_m$ of order 2 with no isolated
intervals.

\medskip

As before we get $O(\log(1/\eps))^{|I_0| + ... + |I_{m-1}|}$ times
an integral nearly identical to what we would have had if starting
with $K_m$. Again there will be extra factors which aid convergence.
The integral in question can be bounded by one of the following
integrals:

\bea \lll{35} \int \int \frac{1}{(1+|p|)^2(1+|q|)^2(1+|p+q|)^3}
e^{-\eps (p^{2}+q^{2})} |p||q|dpdq \\ \nn \int \int
\frac{1}{(1+|p|)^2(1+|q|)^3(1+|p+q|)^2} e^{-\eps (p^{2}+q^{2})}
|p||q|dpdq \\ \nn \int \int \frac{1}{(1+|p|)^3(1+|q|)^2(1+|p+q|)^2}
e^{-\eps (p^{2}+q^{2})} |p||q|dpdq \eea And therefore the following
lemma completes the proof in this case.

\begin{lemma} \lll{bnd6} Each of the integrals in (\ref{35}) is
$o(\log(1/\eps))^2$
\end{lemma}

{\bf Proof:} This is fairly straightforward to calculate using
Lemmas \ref{bnd1} and \ref{bnd4}. For example, by the Cauchy-Schwarz
inequality and symmetry we can bound the first integral by

\be k \lll{yah}\int \int \frac{1}{(1+|p|)^4(1+|p+q|)^3} e^{-\eps
(p^{2}+q^{2})} |p||q|dpdq \ee The $dq$ integral is $(1+|p|)O(1) +
O(\log(1/\eps))$ by Lemma \ref{bnd4}, and thus \ref{yah} is

\be k O(1) \int \int \frac{1}{(1+|p|)^2} e^{-\eps p^2} dp +
O(\log(1/\eps)) \int \int \frac{1}{(1+|p|)^3} e^{-\eps p^2} dp \ee
which is $O(\log(1/\eps))$ by Lemma \ref{bnd1} and the fact that
$\frac{1}{(1+|p|)^3} \in L^1$.

The second and third integrals are identical with $p$ and $q$
interchanged, so we need only do one, let us say the second one.
This is bounded by

\be k \int \int \frac{1}{(1+|p|)(1+|q|)^2(1+|p+q|)^2} e^{-\eps
(p^{2}+q^{2})} dpdq \ee By Cauchy-Schwarz, this is bounded
by

\be \Big(\int \int \frac{1}{(1+|p|)^2(1+|q|)^2} e^{-\eps
(p^{2}+q^{2})}dpdq \Big)^{1/2} \Big( \int \int
\frac{1}{(1+|q|)^2(1+|p+q|)^4} e^{-\eps (p^{2}+q^{2})} dpdq
\Big)^{1/2} \ee This first integral is $O(\log(1/\eps))^2$ by Lemma
\ref{bnd1}, and the second one is $O(\log(1/\eps))$, using Lemma
\ref{bnd1} in conjunction with the fact that

\be \int \frac{1}{(1+|p+q|)^4}dp = O(1) \ee As a simple alternate
proof, one can recall our proof for the case with no isolated
intervals where we constructed the sets $A$ and $B$. Here it is
simple to verify in each case that we can form two sets with the
same properties. The lemma is then proved by the reasoning in steps
(\ref{abs3}) through (\ref{abs6}) above. \qed

\medskip

{\bf Case 3:} There is a $K_m$ consisting of just one interval.

\medskip

Here we must examine in closer detail the proof of Lemma
\ref{remove}. First of all, if there was ever an isolated interval
in some $K_{m'}$ which contained two or more isolated intervals in
$K_{m'-1}$, then the variable corresponding to that interval, say
$u_{(m',j)}$, would have had $l_{m',j} \geq 2$. In that case, by
Lemma \ref{bnd4}, the contribution to (\ref{big7}) of the
$\hat{p}_{m',j}$ integral is $O(1)(1+|u_{(m',j)}-\hat{p}_{m',j})|$.
We see that we can replace the term $O(\log(1/\eps))^{|I_0| + ... +
|I_{m-1}|}$ in (\ref{big8}) with $o(\log(1/\eps))^{|I_0| + ... +
|I_{m-1}|}$, which will finish the proof. Thus we need only consider
the case where $s_1<s_2< ... <s_n < t_n < ... < t_2 < t_1$. In this
case, consider what happens as we remove the first three intervals
(recall that we are assuming that there are at least three
intervals). After removing $(s_n,t_n)$ and then $(s_{n-1},t_{n-1})$
we have

\bea O(\log(1/\eps)) \int e^{-\eps \sum_{1\leq i \leq n-2} p_i^2}
\prod_{1 \leq i \leq n-2} |p_i|
\\ \nn (1+|u_{(2,n-2)}|+O(\log(1/\eps)))\prod_{1 \leq j \leq n_2}
\frac{1}{(1+u_{(2,j)}^2)} \prod_{1 \leq i \leq n-2} dp_i \eea Note
that $u_{2,n-2} = p_1 + ... + p_{n-2}$. We can expand this into two
integrals, namely

\bea \lll{first} O(\log(1/\eps)) \int e^{-\eps \sum_{1\leq i \leq
n-2} p_i^2} \prod_{1 \leq i \leq n-2} |p_i|
\\ \nn (1+|u_{(2,n-2)}|)\prod_{1 \leq j \leq n_2}
\frac{1}{(1+u_{(2,j)}^2)} \prod_{1 \leq i \leq n-2} dp_i \eea and

\bea \lll{sec} O(\log(1/\eps))^2 \int e^{-\eps \sum_{1\leq i \leq
n-2} p_i^2} \prod_{1 \leq i \leq n-2} |p_i|
\\ \nn \prod_{1 \leq j \leq n_2}
\frac{1}{(1+u_{(2,j)}^2)} \prod_{1 \leq i \leq n-2} dp_i \eea The
integral in (\ref{first}) is $O(\log(1/\eps))^{n-2}$ by the same
technique as was used to prove Lemma \ref{remove}. Thus,
(\ref{first}) is $O(\log(1/\eps))^{n-1}$. As for (\ref{sec}), when
we remove the next interval, $(s_{n-2},t_{n-2})$, we have no powers
of $|u_{2,n-2}|$ in the numerator, and by Lemma \ref{bnd4} we do not
pick up an $O(\log(1/\eps))$ term. Thus, (\ref{sec}) is
$o(\log(1/\eps))^n$ as well. This completes the proof of Proposition
\ref{order3}.

\section{Completing the proof}

All that remains is to prove that the processes $\aaa'_{\eps}(T)$
are tight and that the limit process has independent increments.
Both are essentially corollaries of the following lemma:

\begin{lemma} \lll{last} If $b \leq c$, then $(\log(1/\eps))^{-1}\aaa'_{\eps}([a,b]\times[c,d])
\lar 0$ in $L^n$, for any $n \geq 3$.
\end{lemma}

{\bf Proof:} To compute $E[\aaa'_{\eps}([a,b]\times[c,d])]^n$, we
multiply the integrals together as before (see (\ref{expect})). Now,
however, we have $s_i \leq b \leq c \leq t_i$ for all $i$, and it
follows from this that the only configurations of intervals that can
appear here are ones containing just one component of order $n$. We
have shown that these components contribute $o(\log(1/\eps))$ to the
the expectation, and this is enough to prove the lemma \qed

Now that we have this lemma, we can show that the processes
$\aaa'_{\eps}(T)(\log(1/\eps))^{-1}$ are tight. We will show that

\be \lll{tight}
E[(\log(1/\eps))^{-1}(\aaa'_{\eps}(T)-\aaa'_{\eps}(S))]^{2n} \leq
k(T-S)^n \ee where $k$ depends on $n \geq 2$ but can be chosen
independently of $\eps, S$, and $T$, provided they are sufficiently
small. This will prove tightness by, for example, Theorem 12.3 in
\cite{billing}. We can rewrite the left side of (\ref{tight}) as

\be \lll{cot} E[(\log(1/\eps))^{-1}(\aaa'_{\eps}([0,S] \times
[S,T])+\aaa'_{\eps}(D_T \inn \{s,t \geq S \}))]^{2n} \ee We know by
the lemma that $(\log(1/\eps))^{-1}\aaa'_{\eps}([0,S] \times [S,T])
\lar 0$ in $L^{2n}$, so that (\ref{cot}) is bounded by

\be k E[(\log(1/\eps))^{-1}\aaa'_{\eps}(D_T \inn \{s,t \geq S
\})]^{2n} \ee Suppressing the $log$ for the time being, this is
given by

\be \lll{tri} \frac{(-1)^n}{(2\pi)^{4n}} \! \! \int \! \! \!
\int_{D_T^{2n} \inn \{s,t \geq S \}} \! \! \! \! e^{-\eps \sum_j
p_j^2} \prod_{j=1}^n p_{j,1} E\big[ \prod_{j=1}^n e^{i p_j(X_{t_j} -
X_{s_j})}\big] \prod_{j=1}^n ds_j dt_j d^2p_j \ee If we rewrite
$X_{t_j} - X_{s_j}$ as $(X_{t_j} - X_S) - (X_{s_j}-X_S)$, and let
$\beta_t = X_{S+t} - X_S$ be a new Brownian motion this is

\be \lll{tri2} \frac{(-1)^n}{(2\pi)^{4n}} \! \! \int \! \! \!
\int_{D_(T-S)^{2n}} \! \! \! \! e^{-\eps \sum_j p_j^2} \prod_{j=1}^n
p_{j,1} E\big[ \prod_{j=1}^n e^{i p_j(\beta_{t_j} -
\beta_{s_j})}\big] \prod_{j=1}^n ds_j dt_j d^2p_j \ee which is equal
to(reinserting the $log$)

\be k E[(\log(1/\eps))^{-1}\aaa_{\eps}'(T-S)]^{2n} \ee And this is
$O(1) |T-S|^n$, as we showed earlier. This establishes tightness. We
can write

\bea \! \! (\log(1/\eps))^{-1}(\aaa'_{\eps}(T)-\aaa'_{\eps}(S))
\\ \nn = (\log(1/\eps))^{-1}(\aaa'_{\eps}([0,S] \times [S,T])+\aaa'_{\eps}(D_T
\inn \{s,t \geq S \})) \eea Since
$(\log(1/\eps))^{-1}\aaa'_{\eps}([0,S] \times [S,T]) \lar 0$ and
$\aaa'_{\eps}(D_T \inn \{ s,t \geq S \}) \in \unn \{ \sigma
(X_t-X_s) : S \leq s,t \leq T \}$, we see that $(\log(1/\eps))^{-1}
\aaa'_{\eps}(T)$ has asymptotically independent increments. This
shows that the limit process, $W_T$, has independent increments, and
completes the proof of Theorem 1.

\section{Symmetric stable processes}

We will now prove Theorem \ref{stab}. The proof of this theorem is,
naturally, very similar to the proof in the Brownian motion case, so
we will in many cases just refer to steps undertaken in the previous
proof. In particular, the general outline (Sections 2 and 5) is
identical in both cases; the only difference lies in some of the
calculations.

The main difficulty is in showing that the integrals corresponding
to components of order two converge. $X_t$ is a symmetric stable process of index $\bb$ where $1 < \bb < 2$. The density of
$X_t$ is given by

\be \lll{psb2} f_{\eps}(x) = \frac{1}{(2\pi)^{2}} \int e^{ipx-\eps p^{\bb}}
d^{2}p \lll{2}\ee Thus,

\be \lll{psb3} f_{\eps}'(x) = \frac{i}{(2\pi)^{2}} \int p_1 e^{ipx-\eps p^{2}}
d^{2}p \lll{2}\ee Proceeding as in section 3, the
first integral is

\bea \lll{11a} \int \int
\frac{(1-e^{-p^{\bb}})}{p^{\bb}}\frac{(1-e^{-q^{\bb}})}{q^{\bb}}
\frac{(1-e^{-(p+q)^{\bb}})}{(p+q)^{\bb}} e^{-\eps(p^{\bb}+q^{\bb})}
p_{1}q_{1}dpdq \\ = \eps^{3-6/\bb}\int \int \frac{(1-e^{-p^{\bb}
/\eps})}{p^{\bb}}\frac{(1-e^{-q^{\bb} /\eps})}{q^{\bb}}
\frac{(1-e^{-(p+q)^{\bb} /\eps})}{(p+q)^{\bb}}
e^{-(p^{\bb}+q^{\bb})} p_{1}q_{1}dpdq \nn \eea In order to prove
that this integral converges as $\eps \lar 0$, it is enough to show
that

\be \int \int \frac{1}{p^{\bb-1}}\frac{1}{q^{\bb-1}}
\frac{1}{(p+q)^{\bb}} e^{-(p^{\bb}+q^{\bb})}dpdq \ee converges, and
then to apply the dominated convergence theorem. We need only
consider the integral over $\{|p|, |q| < 1\}$, for in order to
evaluate the integral over, say, $A=\{|p|>1\}$ we may divide $A$
into the disjoint union of $B= \{|q|<1/2\} \inn A$, $C =
\{|p+q|<1/2\} \inn A$, and $D = A- (B \unn C)$. The integrals over
$B$ and $C$ are both bounded because $\bb < 2$ i.e. we have only
integrable singularities. And the integral over $D$ is bounded by a
constant times

\be \int \int e^{-(p^{\bb}+q^{\bb})}dpdq < \infty \ee So we must
consider the integral

\be \lll{sun} \int \int_{\{|p|,|q|<1\}}
\frac{1}{p^{\bb-1}}\frac{1}{q^{\bb-1}} \frac{1}{(p+q)^{\bb}}dpdq \ee
We manipulate the integral as follows:

\bea \int _{|p|<1} \frac{1}{p^{\bb-1}} \int_{\{|q|<1\}}
\frac{1}{q^{\bb-1}} \frac{1}{(p+q)^{\bb}}dqdp
\\ \nn = \int_{|p|<1} \frac{1}{p^{3 \bb -2}} \int_{\{|q|<1\}}
\frac{1}{(q/|p|)^{\bb-1}} \frac{1}{(p/|p| + q/|p|)^{\bb}} dqdp \eea
The argument of $p$ (thought of as a complex number) is irrelevant,
so we may replace $p/|p|$ by $1$, and substitute $q' = q/|p|$ to get

\be \lll{67} \int_{|p|<1} \frac{1}{p^{3 \bb -4}} \Big(
\int_{\{|q|<1/|p|\}} \frac{1}{q^{\bb-1}} \frac{1}{(1 + q)^{\bb}} dq
\Big) dp \ee If $\bb > 3/2$ then the $dq$ integral is bounded
independently of $|p|$ (since then $\frac{1}{q^{\bb-1}} \frac{1}{(1
+ q)^{\bb}} \in L^1$), so that (\ref{sun}) is bounded by a constant
times

\be \int_{\{|p|<1\}} \frac{1}{p^{3\bb - 4}} dp \ee which is finite,
as $\bb < 2$. If $\bb < 3/2$ (resp. $\bb = 3/2$), then the $dq$
integral in (\ref{67}) is $O(|p|)^{2\bb-3}$ (resp. $O(|\log|p||))$,
so that (\ref{sun}) is bounded by a constant times

\be \int_{\{|p|<1\}} \frac{1}{p^{\bb - 1}} dp \ee when $\bb < 3/2$
and

\be \int_{\{|p|<1\}} \frac{|\log(|p|)|}{p^{1/2}} dp \ee when $\bb =
3/2$. These integrals are both finite.

The second configuration of intervals gives rise to the following:

\be \lll{fil} \int \int \frac{(1-e^{-p^{\bb} /\eps})^2}{p^{2
\bb}}\frac{(1-e^{-(p+q)^{\bb} /\eps})}{(p+q)^{\bb}}
e^{-(p^{\bb}+q^{\bb})} p_{1}q_{1}dpdq \ee This integral is more
difficult as for some $\bb$ the integrand is not in $L^1$ were we to
remove the terms involving $\eps$ (there is a non-integrable
singularity at $p=0$ when $\bb \geq 3/2$). We will first show that
(\ref{fil}) is bounded independently of $\eps$. We isolate the $dq$
integral:

\be \lll{mm} \int \frac{(1-e^{-(p+q)^{\bb} /\eps})}{(p+q)^{\bb}}
e^{-q^{\bb}} q_{1} dq \ee We will show that this is $|p|O(1)$ (the
$O$ here refers to $\eps$). Because we will refer to this result
later, we isolate it as a lemma (which we state in slightly greater
generality).

\begin{lemma} \lll{bnd11} For any $a$ with $0<a<T$ and any $p \in {\bf R}^2$, we have

\be \lll{ozo} \Big|\int
\frac{1-e^{-(p+q)^{\bb}a/\eps}}{(p+q)^{\bb}}e^{q^{\bb}} q_1 dq\Big|
= |p|O(1) \ee independently of $a$.

\end{lemma}

Proof: We can drop the $e^{-(p+q)^{\bb}a/\eps}$ term. (\ref{ozo}) is
bounded by

\be \lll{mm1} \Big| \int \frac{1}{q^{\bb}} e^{-(q-p)^{\bb}}
(q_{1}-p_1) dq \Big| \ee Expand the $(q_1 - p_1)$ term. The second
term is bounded by

\be \lll{bg}|p| \int \frac{1}{q^{\bb}}e^{-(q-p)^{\bb}} dq \ee The
integrand is bounded by the function

\be \frac{1}{q^{\bb}}1_{\{|q|<1\} } + e^{-(q-p)^{\bb}} 1_{\{|q|\geq
1\} } \ee which is bounded in $L^1$ independently of $p$. Thus,
(\ref{bg}) is $|p|O(1)$. To bound the first term we subtract

\be \int \frac{1}{q^{\bb}} e^{-q^{\bb}}q_1 dq \ee which is $0$ by
symmetry. This gives us

\bea \lll{ccc} \Big| \int \frac{1}{q^{\bb}}
(e^{-(q-p)^{\bb}}-e^{-q^{\bb}}) q_{1} dq \Big| \\ \nn \leq \int
\frac{1}{q^{\bb}} |e^{-(q-p)^{\bb}}-e^{-q^{\bb}}| |q| dq \eea We
split this up into the integral over the region $\{|q|<2|p|\}$ and
$\{|q| \geq 2|p|\}$. The integral over the first region is bounded
by

\bea k \int_{\{|q|<2|p|\} } |q|^{1-\bb} dq \\ \nn = k
\int_0^{2|p|}r^{2-\bb} dr \\ \nn = k |p|^{3-\bb} \eea Here $k$ is a
constant which may change from line to line. This is
$(|p|+|p|^2)O(1)$. On the region $\{|q| \geq 2|p|\}$ suppose first
that $e^{-(q-p)^{\bb}} \geq e^{-q^{\bb}}$. Then

\bea |e^{-(q-p)^{\bb}}-e^{-q^{\bb}}| \\ \nn \leq
e^{-(|q|-|p|)^{\bb}}-e^{-q^{\bb}} \\ \nn = \bb \int_{|q|-|p|}^{|q|}
x^{\bb -1}e^{-x^{\bb}}dx \\ \nn \leq
k|p||q|^{\bb-1}e^{-(|q|-|p|)^{\bb}}\eea The last inequality is the
length of the interval being integrated over multiplied by a term
which bounds the integrand. Plugging this into (\ref{ccc}) gives a
bound of

\be k|p| \int_{|q| \geq 2|p|} e^{-(|q|-|p|)^{\bb}} dq \leq k|p|
\int_{|q| \geq 2|p|} e^{-(|q|/2)^{\bb}}dq = O(|p|) \ee In the case
$e^{-(q-p)^{\bb}} < e^{-q^{\bb}}$ we have

\bea |e^{-(q-p)^{\bb}}-e^{-q^{\bb}}| \\ \nn \leq
e^{-|q|^{\bb}}-e^{-(|q|+|p|)^{\bb}} \\ \nn = \bb
\int_{|q|}^{|q|+|p|} x^{\bb -1}e^{-x^{\bb}}dx \\ \nn \leq
k|p|(|q|+|p|)^{\bb-1}e^{-|q|^{\bb}}\eea Since $|q| \geq 2|p|$ this
is $k|p||q|^{\bb-1}e^{-|q|^{\bb}}$ Thus, the contribution to
(\ref{ccc}) of this region is bounded by

\be k|p| \int e^{-|q|^{\bb}} dq  = |p|O(1) \ee This shows that

\be \lll{mmn} \int \frac{1}{q^{\bb}} e^{-(q-p)^{\bb}} q_{1} dq =
(|p| + |p|^2)O(1)\ee It is also $O(1)$, however, since the integrand
is bounded by

\be \frac{1}{q^{\bb-1}}1_{\{|q|<1\}} + e^{-(q-p)^{\bb}} \ee which is
bounded in $L^1$ independently of $p$. So (\ref{mmn}) is $|p|O(1)$
for $p$ small, and $O(1)$ for $p$ large. We conclude that
(\ref{mmn}) is $|p|O(1)$ for all $p$. \qed

\medskip

This lemma allows us to see that (\ref{fil}) is bounded by

\be \lll{off} k \int \frac{(1-e^{-p^{\bb} /\eps})^2}{p^{2 \bb}}
e^{-p^{\bb}}|p|^2 dp \ee The extra powers of $p$ in the numerator
are enough to convert our singularity at $0$ into an integrable one,
and it follows that (\ref{off}) is bounded by

\be \lll{ff} k \int \frac{1}{p^{2 \bb}} e^{-p^{\bb}}|p|^2 dp <
\infty \ee We have showed that (\ref{fil}) is bounded independently
of $\eps$. This alone does not show that (\ref{fil}) converges.
However, convergence is proved using the same ideas, as follows. Let
the value of (\ref{fil}) be denoted by $A(\eps)$. We will show that,
for any $\dd > 0$, there is an $\eps '>0$ such that if $0< \eps_1,
\eps_2 < \eps '$ then $|A(\eps_1) - A(\eps_2)|<\dd$. This will prove
convergence. We will assume below that $0<\eps_1<\eps_2<\eps'$. We
have

\bea \lll{fil2} A(\eps_1) - A(\eps_2) \\ = \nn \int \int \Big{(}
\frac{(1-e^{-p^{\bb} /\eps_1})^2}{p^{2 \bb}}\frac{(1-e^{-(p+q)^{\bb}
/\eps_1})}{(p+q)^{\bb}} - \frac{(1-e^{-p^{\bb} /\eps_2})^2}{p^{2
\bb}}\frac{(1-e^{-(p+q)^{\bb} /\eps_2})}{(p+q)^{\bb}}\Big{)} \\ \nn
e^{-(p^{\bb}+q^{\bb})} p_{1}q_{1}dpdq \eea We will rewrite the
difference

\be (1-e^{-p^{\bb} /\eps_1})^2 (1-e^{-(p+q)^{\bb} /\eps_1}) -
(1-e^{-p^{\bb} /\eps_2})^2 (1-e^{-(p+q)^{\bb} /\eps_2}) \ee as

\bea (1-e^{-p^{\bb} /\eps_1})^2 [(1-e^{-(p+q)^{\bb} /\eps_1}) -
(1-e^{-(p+q)^{\bb} /\eps_2})] \\ \nn + [(1-e^{-p^{\bb} /\eps_1})^2 -
(1-e^{-p^{\bb} /\eps_2})^2] (1-e^{-(p+q)^{\bb} /\eps_2}) \eea and
handle each term in this sum separately. The first one gives rise to
the integral

\be \lll{fil2} \int \int \frac{(1-e^{-p^{\bb} /\eps_1})^2}{p^{2
\bb}}\frac{(e^{-(p+q)^{\bb} /\eps_2}-e^{-(p+q)^{\bb}
/\eps_1})}{(p+q)^{\bb}} e^{-(p^{\bb}+q^{\bb})} p_{1}q_{1}dpdq \ee As
in step (\ref{mm1}) the $dq$ integral is

\be \lll{mm2} \int \frac{e^{-q^{\bb} /\eps_2}(1-e^{-q^{\bb}
/\eps_3})}{q^{\bb}} e^{-(q-p)^{\bb}} (q_{1}-p_1) dq \ee where
$\eps_3 = \frac{\eps_1 \eps_2}{\eps_2 - \eps_1}>0$ This is in turn
bounded by

\be \lll{mm3} \int e^{-q^{\bb}
/\eps'}\frac{1-e^{-q^{\bb}/\eps_3}}{q^{\bb}} e^{-(q-p)^{\bb}}
(q_{1}-p_1) dq \ee We may now follow steps (\ref{mm1}) through
(\ref{ff}), and it is straightforward to verify in each case that
the extra $e^{-q^{\bb} /\eps'}$ term allows us to replace the $O(1)$
by $o(1)$ (the $o$ now refers to $\eps'$). This implies that
(\ref{fil2}) can be made arbitrarily small by choosing $\eps'$
sufficiently small. As for the second integral

\be \lll{fil3} \int \! \! \! \int \! \! \frac{((1-e^{-p^{\bb}
/\eps_1})^2 \! \! - \! (1 \! \! -\! \! e^{-p^{\bb}
/\eps_2})^2)}{p^{2 \bb}} e^{-(p^{\bb}+q^{\bb})} p_{1}q_{1}dpdq \ee
We can rewrite $((1-e^{-p^{\bb} /\eps_1})^2 - (1 - e^{-p^{\bb}
/\eps_2})^2)$ as

\be ((1-e^{-p^{\bb} /\eps_1}) + (1 - e^{-p^{\bb}
/\eps_2}))e^{-p^{\bb}/\eps_2}(1-e^{-p^{\bb}/\eps_3}) \ee and we see
that we can bound (\ref{fil3}) by

\be k \int \frac{e^{-p^{\bb}/\eps'}|p|e^{-p^{\bb}}}{p^{2\bb}} \Big|
\int \frac{(1-e^{-(p+q)^{\bb} /\eps_2})}{(p+q)^{\bb}}e^{-q^{\bb}}q_1
dq \Big| dp \ee We have shown above that the $dq$ integral is $(|p|
+ |p|^2)O(1)$, and that this implies that the entire integral
converges. Furthermore, as $\eps' \lar 0$, the dominated convergence
theorem implies that the value of the integral approaches zero.
Again we see that if we choose $\eps'$ sufficiently small we can
make (\ref{fil3}) arbitrarily small. This shows that if $\eps_n$ is
a sequence converging to zero then $A(\eps_n)$ converges. Thus,
$\lim_{\eps \lar 0} A(\eps)$ exists, and we define

\be \int \int \frac{1}{p^{2\bb}} \frac{1}{(p+q)^{\bb}}
e^{-(p^{\bb}+q^{\bb})}p_1 q_1 dpdq \ee to be this limit. This
completes the calculation for components of order 2.

For a component of order $n \geq 3$ we have the following integral:

\be \lll{biga} \int e^{-\eps \sum p_i^{\bb}} \prod_i (p_i)_1 \Big(
\int_{\sum c_j < T} \prod_j e^{-u_j^{\bb} c_j} \prod_j dc_j \Big)
\prod_i dp_i \ee We must show that this is $o(\eps^{-(3n/\bb -
3n/2)})$. This would be a bit of a chore were it not that we have
done almost all of the work already in the Brownian motion case. For
instance, suppose we have a configuration with no isolated
intervals. Then (\ref{biga}) can be bounded by(see (\ref{abs}))

\bea \lll{absa} \int \frac{e^{-\eps \sum p_i^{\bb}} \prod
|p_i|}{\prod_{j=1}^{2n-1}(1+|u_j|)^{\bb}}\prod dp_i \\ \nn =
\eps^{-(3n/\bb -2n - 1)} \int \frac{e^{-\sum p_i^{\bb}} \prod
|p_i|}{\prod_{j=1}^{2n-1}(\eps+|u_j|)^{\bb}}\prod dp_i \eea We are
done if we can bound this integral effectively. We know from earlier
work that if $\bb$ were replaced by 2 in this integral then it would
be $O(\log(1/\eps))^n$, which is certainly good enough. We can bound
as follows using Holder's inequality:

\bea \int \lll{supe} \frac{e^{-\sum p_i^{\bb}} \prod
|p_i|}{\prod_{j=1}^{2n-1}(\eps+|u_j|)^{\bb}}\prod dp_i  \\ \nn \leq
\Big( \int \frac{e^{-\sum p_i^{\bb}} \prod
|p_i|}{\prod_{j=1}^{2n-1}(\eps+|u_j|)^{2}}\prod dp_i \Big)^{\bb/2}
\Big( \int e^{-\sum p_i^{\bb}} \prod |p_i|\prod dp_i
\Big)^{(2-\bb)/2} \eea A quick examination of the proofs of Lemmas
\ref{bnd1}, \ref{bnd2}, \ref{bnd4}, and \ref{bnd6} will show that
the conclusions of these lemmas remain valid if any $e^{-p^2}$'s in
the hypotheses are replaced by $e^{-p^{\bb}}$. We can conclude that
(\ref{supe}) is $O(\log(1/\eps))^n$, and this component is therefore
sufficiently bounded. We do the same thing in the isolated interval
case, with Lemma \ref{bnd3} replaced by Lemma \ref{bnd11}. This
completes the proof of Theorem \ref{stab}.

\medskip

\noindent{\bf Acknowledgements}

I am deeply indebted to my advisor Jay Rosen, who suggested this
problem to me, who taught me a great deal, and whose help and
generosity were invaluable in completing this work.


\def\noopsort#1{} \def\printfirst#1#2{#1} \def\singleletter#1{#1}
   \def\switchargs#1#2{#2#1} \def\bibsameauth{\leavevmode\vrule height .1ex
   depth 0pt width 2.3em\relax\,}
\makeatletter \renewcommand{\@biblabel}[1]{\hfill#1.}\makeatother

\end{document}